\documentclass[journal,12pt,onecolumn]{IEEEtran} 
\usepackage[english]{babel}
\usepackage{csquotes}
\usepackage[backend=biber,style=authoryear,uniquename=init,giveninits]{biblatex}
\usepackage{hyperref}  
\usepackage{subfiles}
\usepackage{longtable}
\usepackage{booktabs,subcaption,amsfonts,dcolumn}

\newcounter{eqc} 

\newif\ifsubbuild 
\subbuildtrue 

\newif\ifparaen
\ifparaen 
    \newcommand{\keyq}[1]{\noindent \color{gray} \textit{Key question(s): #1} \color{black} }
    \newcounter{pitemcntr} 
    \newcommand{\pitem}[1]{\refstepcounter{pitemcntr} \par \textbf{\thepitemcntr. #1: }} 
    \newcommand{\pitemo}[2]{\refstepcounter{pitemcntr} \par \textbf{\thepitemcntr. #1:} #2} 
    \newenvironment{parag}[1]
        { \keyq{#1} }
        { \par \setcounter{pitemcntr}{0} 
        \vspace{0.6cm}
        }
\else 
    \newdimen\origiwspc%
    \origiwspc=\fontdimen2\font
    \newcommand{\keyq}[1]{}
    \newcommand{\pitem}[1]{\fontdimen2\font=\origiwspc}
    \newcommand{\pitemo}[2]{\fontdimen2\font=\origiwspc{}#2\fontdimen2\font=-0.01ex}
    \newenvironment{parag}[1]{\fontdimen2\font=-\origiwspc}{\fontdimen2\font=\origiwspc} 
\fi

\hypersetup{
    colorlinks=true,
    linkcolor=blue,
    filecolor=magenta,      
    urlcolor=cyan,
}

\usepackage{amsmath} 
\usepackage{amsfonts}
\usepackage{calrsfs} 

\usepackage{blkarray} 
\usepackage{ragged2e} 
\usepackage{enumitem} 
\usepackage{setspace}
\usepackage{xcolor} 
\usepackage{color,soulutf8}
\usepackage{graphicx} 
\usepackage{caption}
\usepackage{subcaption}
\usepackage{bbm}
\usepackage[linesnumbered,ruled]{algorithm2e} 
\usepackage{algorithmic}
\usepackage{lscape}
\usepackage[T1]{fontenc}
\usepackage[utf8]{inputenc}
\usepackage[normalem]{ulem}
\usepackage{pdfpages}
\usepackage{booktabs, multirow} 
\usepackage{soul}
\usepackage{changepage,threeparttable} 
\usepackage{mathtools}

\usepackage{imakeidx}
\makeindex[columns=3, title=Alphabetical Index, options= -s example_style.ist, intoc]

\usepackage{pifont}

\newcommand{\indx}[1]{\hypertarget{#1}{\textit{#1\index{#1}}}} 
\newcommand{\indxy}[2]{\hypertarget{#2}{\textit{#2\index{#1#2}}}} 
\newcommand{\rindx}[1]{\hyperlink{#1}{\color{black!80!blue}{\textit{#1}}}}
\newcommand{\rindxy}[2]{\hyperlink{#1}{\color{black!80!blue}{\textit{#2}}}}

\newtheorem{thm}{Theorem}[section]
\newtheorem{lem}[thm]{Lemma}
\newtheorem{corr}{Corollary}[thm]
\newtheorem{rem}{Remark}

\newenvironment{defy}[1]{\indx{#1}: }{}
\newenvironment{defxy}[2]{\indxy{#1}{#2}: }{}
\newenvironment{prf}{\noindent\textit{Proof. }}{\prvd}
\newcommand{\prvd}{\ding{114}}

\newcommand{\autha}[1][*]{\textcolor{black}{#1}} 
\newcommand{\authar}[1][*]{\textcolor{black}{#1}} 

\newcommand{\hgl}[1]{#1} 

\newcommand{\work}[1]{work}

\subbuildfalse 

\begin{document}

\title{Temporal Recurring Unavailabilities in Multi-agent Rural Postman Problem: Navigating railway tracks during availability time intervals} 
\author{Somnath Buriuly, Leena Vachhani, Arpita Sinha, Sivapragasam Ravitharan, Sunita Chauhan
\thanks{S. Buriuly is a postdoctoral fellow with CoE-OGE, IIT Bombay, India}
\thanks{L. Vachhani and A. Sinha are with Systems and Control Engineering,
IIT Bombay, India}
\thanks{S. Ravitharan is with the Department of Mechanical and Aerospace Engineering, Monash University, Australia}
\thanks{S. Chauhan is with the Center for Equitable $\&$ Personalized Health, Plaksha University, India}}






\markboth{Manuscript for submission}%
{Buriuly \MakeLowercase{\textit{et al.}}: Temporal Recurring Unavailabilities in Multi-agent Rural Postman Problem: Navigating railway tracks during availability time intervals} 

\maketitle

\begin{abstract}
	Time-dependent (or temporal) properties may arise in many network-based planning problems, particularly in the routing and scheduling of railway track inspection problems. The availability of tracks depends on the train schedules, maintenance possessions, etc. 
	In the absence of side constraints, this routing and scheduling problem is formulated as a multi-agent rural postman problem on a temporal-directed network; where a given set of rail track sections must be visited while respecting the temporal attributes due to railway track unavailabilities. 
	In this work, we adopt a three-index formulation for the multi-agent Rural Postman Problem with Temporal Recurring Unavailabilities (RPP-TRU) and frame it as a Mixed Integer Linear Programming (MILP) problem. In addition, we propose relevant theoretical studies for RPP-TRU to ensure the feasibility of the proposed optimization problem. Two approaches of an exact algorithm are proposed, based on Benders' decomposition framework, to address the disjunctive unavailability constraints occurring in its scheduling sub-problems, alongside the NP-Hard routing (master) problem. A polynomial-time algorithm is designed to address the scheduling sub-problem, while the NP-Hard master problem is solved using MILP toolbox. 
	Comparison results with RPP (without temporal constraints) show a minor compromise with the spatial cost solution with significantly less delay, hence suitable for real-world routing and scheduling applications occurring in a shared network like railways. 
	A simulation study on a part of the Mumbai suburban railway network demonstrates the working of the proposed methodology under a realistic setting. 
\end{abstract}


\begin{parag}{\textbf{Asked ChatGPT - What are the key points of the abstract, listed in brief?}}
    \pitemo{The problem of routing and scheduling in railway track inspection involves time-dependent constraints due to track unavailabilities.}{}
    \pitemo{The problem is formulated as a multi-agent Rural Postman Problem with Temporal Recurring Unavailabilities (RPP-TRU) in a Mixed Integer Linear Programming (MILP) framework.}{}
    \pitemo{Two approaches based on Benders' decomposition framework are proposed to solve the scheduling sub-problems and the NP-Hard routing problem.}{}
    \pitemo{The proposed methodology is shown to have a minor trade-off on spatial cost solution but significantly less delay compared to RPP without temporal constraints.}{}
    \pitemo{The methodology is demonstrated through a simulation study on a part of the Mumbai suburban railway network.}{}
\end{parag}

\begin{parag}{\textbf{Questions generated by ChatGPT and Me, for this paragraph -} What is the motivation of this work? What is the problem being addressed in this work? Why is studying RPP-TRU important for railway networks with repeating train schedules?} 
    \pitemo{Motivated from planning on shared infrastructure}{Urbanization has increased the demand for railway services thereby increasing wear and tear of the infrastructural resources. To satisfy the transportation demand while maintaining infrastructure reliability, frequent monitoring is necessary, while ensuring least disruption to the regular train schedules. } 
    \pitemo{This work models such planning problems}{In this work, we address inspection planning of railway tracks on a time-shared environment. The core problem is a \autha[Rural Postman Problem] (RPP) with time-dependencies occurring due to sharing of railway infrastructure for various services. The RPP is an NP-Hard problem (see \cite{corberan2015}) that - in context of our \work{} - translates into the routing of inspection wagons for servicing a subset of the track section, while the time-dependencies on the arcs of RPP represent the unavailability of tracks caused by train schedules. }
    \pitemo{Practical significance is observed due to recurring train schedules}{In addition, recurrence is an anticipated property in railway networks with repeating train schedules. Hence, we study a variant of RPP coined as Rural Postman Problem with Temporal Recurring Unavailabilities (RPP-TRU). }
\end{parag} 

\begin{parag}{What exists in the literature?} 
    \pitemo{Few factors require study of temporal attributes}{Temporal (or time-dependent) attributes have been studied in the research literature for many network-based problems. Some of the factors that lead to temporal problems are congestion in roadways, customer time windows, unavailabilities due to rail schedule, etc. } 
    \pitemo{There is literature for roadway temporal problems}{\cite{cordeau} and \cite{calogiuri} present \autha[Time-Dependent Traveling Salesman Problem] (TDTSP) and Time-Dependent Rural Postman Problem (TDRPP) respectively, that reflects temporal attributes in roadway network due to traffic congestion. The temporal attributes of railway networks are contrasting in comparison to the roadway problems. In roadway networks, the travel time is affected by traffic congestion, while in railways the connectivity of the underlying network is lost temporarily. }
    \pitemo{The literature for railway temporal problems is different}{\cite{lannez} presents a routing and scheduling problem for the inspection of railway tracks using ultrasonic sensor-equipped inspection vehicles, where the temporal attributes are modeled with integer variables using discretized unavailability calendars. The unavailability window restrictions are a generalization of the time window constraints, as unavailabilities may occur multiple times and over arbitrary arcs; while in problems with time window constraints, the restrictions are imposed for a single interval only. \cite{bogue2022} proposes a generalization of a routing problem with multiple time windows, however, the windows affect the availability of the service vertices/arcs only. \cite{buriuly2022} proposes a three-index formulation for CARP-TU and RPP-TU, that models the railway routing problem using integer variables, and the block unavailability using continuous variables i.e. without discretization. }
    \pitemo{We have adopted the formulation from the literature}{In this work, we have taken the three-index formulation from RPP-TU, and applied it to the RPP-TRU. The benefit of this approach is achieved by overcoming the dependency of the problem size (number of variables and constraints) on time discretization and planning horizon. The approach is suitable for busy networks that have fast trains running on smaller length signal blocks (track section between two signals) with frequent signal block unavailabilities. }
\end{parag}


The RPP-TRU formulation utilizes temporal variables to address the unavailability of track sections for inspection. Unlike existing work, we show 
a stronger connectivity condition, a theoretical guarantee for the existence of a solution, as well as data pruning  to merge (thereby reduce) the number of unavailability intervals. We also propose two variants of exact algorithms, built on Benders' decomposition framework, to improve over simple MILP-based solvers 
with the aim for time improvements over sub-optimal solutions that ignore temporal attributes of the RPP-TRU. We limit our study 
for medium-sized sparse networks (up to 55 vertices with 66 arcs and 20 tasks) because the computation time for an optimal solution of NP-Hard problems is exponentially proportional to the size of the network. 

\begin{parag}{\textbf{Asked ChatGPT to list the contributions in the previous paragraph} List of contributions:}
 \pitemo{Adapting the RPP-TU framework from \cite{buriuly2022} for the proposed RPP-TRU (RPP-TU with recurring unavailabilities).}{}
 \pitemo{Utilizing temporal variables to address the unavailability of track sections in the RPP-TRU formulation.}{}
 \pitemo{Providing stronger connectivity conditions for the theoretical guarantee of a solution, as well as data cleaning to reduce the number of unavailability intervals.}{}
 \pitemo{Proposing two variants of exact algorithms, built on Benders' decomposition framework, to improve over simple MILP-based solvers.}{}
 \pitemo{Showing time improvements over sub-optimal solutions that ignore temporal attributes of the RPP-TRU, with a minimal trade-off in fuel/spatial cost, for medium-sized sparse networks (up to 55 vertices with 66 arcs and 20 tasks).}{}
\end{parag}

The major contributions of this work are summarized as follows: 
\begin{itemize}
    \item 
    A theoretical analysis to guarantee the existence of a solution in the multi-agent Rural Postman Problem with Temporal Recurring Unavailabilities (RPP-TRU). 
    \item Attributes like well-defined temporal dataset and recurrence in the unavailability of arcs are investigated from the application perspective of RPP-TRU, and supported through a simulation study on Kurla-Thane-Vashi-Kurla (Mumbai) suburban railway network. 
    \item Two approaches of the proposed algorithm are developed for RPP-TRU, based on Benders' decomposition algorithm framework. This work demonstrates the non-trivial equivalence between the proposed cuts for RPP-TRU and the well-known form of Benders' cuts from the literature; thereby ensuring convergence to an optimal solution given sufficient computational time. 
\end{itemize}

We survey relevant works from the literature in Section \ref{rpptu:sec:lit}. Preliminaries for RPP-TRU is discussed in Section \ref{rpptu:sec:prelims}. The three-index formulation for RPP-TRU is presented in Section \ref{rpptu:sec:formulate} along with the theoretical analysis of temporal recurrence. In Section \ref{rpptu:sec:algo} two approaches of an exact algorithm are proposed, for obtaining an optimal solution of RPP-TRU. The methodology is supported by simulation studies and comparison studies presented in Section \ref{rpptu:sec:resdisc}. Lastly, the work is concluded with a summary of the main findings and its projection into future research. 

\section{Literature survey }
\label{rpptu:sec:lit}

\ifsubbuild
\setcounter{section}{4}
\section{Literature survey}
\fi

\begin{parag}{What does the routing and scheduling problem comprise of? How are various routing and scheduling problems involving railways and roadways tackled in the literature? }
 \pitemo{Railway and roadway network planning problems are sub-categorized under routing and scheduling problems, and they are modeled as constrained optimization problems}{
 Planning problems over a railway or roadway network often fall in the category of routing and/or scheduling problems. These planning problems are optimization problems that minimize costs like fuel consumption, time taken, etc. while ensuring that the numerical solution of the optimization problem is valid for the given network. Hence, a suitable set of constraints - e.g. flow constraints, service constraints, etc. - is introduced into the optimization problem to achieve practically feasible routes which is a valid plan for traversing the network.
 }
 \pitemo{There are a few real-world railway network planning problems studied in the literature}{
 Many practical problems arising in railways involve routing and/or scheduling of trains or wagons in the railway network. Train time-tabling is one such problem where the route and schedule of trains need to be determined in an optimal way while respecting network restrictions and safety regulations; see \cite{peeters, thomas}. Inspection and maintenance activities for railway infrastructure also require routing and scheduling of agents on the network; see \cite{lannez, peng, pengThesis, liden, budai}.
 }
 \pitemo{Unconventional - we resepct the unavailability time interval}{Usually, the inspection services are only conducted at night when the train services are nonoperational, see \cite{chen2023}. In this work, we explore the application where the inspection may be conducted during the operating hours of the regular trains.}
\end{parag}

\begin{parag}{How is management done in railways? What are the routing and scheduling practices in railway management?}
 \pitemo{There are many challenges in railway network planning problems which insist for distributed management}{
 Railways feature a vast infrastructure capable of serving multiple transportation objectives. Capacity restrictions of the infrastructure and increasing demand for transportation services motivate the sharing of this resource. This leads to multiple inter-connected transportation challenges requiring the transportation of people and goods, preventive maintenance planning, etc. Some of these challenges are presented as optimization problems, like efficient utilization of allocated resources and optimally sharing them among different units. Efficient management of this vast railway network is extremely difficult for any isolated computational tool. Hence, a distributed approach is the most preferred zeroed-choice, giving rise to hierarchical optimization problems where important ones are solved first.
 }
 \pitemo{Infrastructure Manager supervises the distributed measurement}{The Infrastructure Manager (IM) decides the hierarchy, segregates sub-problems, and assigns them to suitable servicing units called Railway Undertakings (RUs). \cite{thomas} and \cite{peeters} report various planning phases in detail while introducing the railway system.
 } 
 \pitemo{Takeaway}{The routing and scheduling of agents for railway track health monitoring is generally handled by an RU, by requesting possession of railway track signal blocks from IM.
 } 
\end{parag}
 
\begin{parag}{How is the mathematical modeling of railway network planning problems addressed in the literature?} 
 \pitemo{Common practices for modeling these problems in literature involved discrete optimization}{
 Mathematical modeling of these routing and scheduling problems involve representing the discrete decisions/choices of the agents as integer/binary variable and then incorporating the practical restrictions as constraints expressed in terms of the discrete decisions.
 }
 \pitemo{These discrete optimization problems are further sub-categorized into VRPs and ARPs}{
 Almost all routing and scheduling problems are categorized as Vehicle Routing Problems (VRPs) or Arc Routing Problems (ARPs). 
 The main focus in VRPs and ARPs are point locations (vertices) and inter-relations linking these points (arcs), respectively. The primary objective in these routing problems is to optimize the traversal cost of the vehicle(s) such that a subset of either vertices (for VRP) or arcs (for ARP) of the network (graph) are visited at least once. The general problem is very ideal and ignores many practical restrictions like network and vehicle capacity, inspection and maintenance time slots, heterogeneous vehicle characteristics, etc. }
 \pitemo{Many variants of ARP and VRP exist for modeling the railway network problems - they are all NP-Hard}{Based on some of the practical considerations, many variants of VRPs and ARPs have been introduced in the literature, like VRP with Time-Window (TW), Capacitated Arc Routing Problem (CARP), etc. discussed in \cite{aeiseltLaporte}, \cite{beiseltLaporte}, \cite{toth2014}, \cite{corberan2015}, etc. Both the general problem (mostly NP-Hard) as well as the practical problem (with many side constraints) - proposing exact solutions as well as heuristics - have been attempted. 
 }
\end{parag} 


\keyq{What research is available for time-dependent problems in railways and roadways? }
Routing and scheduling problems in railways and roadways have shared infrastructure, leading to time-varying attributes in the VRP or ARP formulation. 
A challenging problem with time-varying attributes has been studied by \cite{calogiuri}, which solves a Time-Dependent Rural Postman Problem (TDRPP) in a road traffic environment. 
TDRPP for roadway networks uses the temporal variant of the graph that is traversable if the underlying graph is connected. However, railway-based unavailability scenarios can result in disconnected sub-graphs at various time instants. For instance, a track section, which is available for inspection between $0$ and $T$ time units only, can never be traversed by a slow inspection wagon (agent) that has a traversal time of more than $T$ time units for this track section. Hence a stronger connectivity condition for agent-network compatibility is required. Furthermore, TDRPP addresses a single-agent scenario and considers a {First-In-First-Out} (FIFO) rule wherein if a vehicle leaves a vertex (say $u$) for another vertex (say $v$) at a given time, leaving vertex $u$ at a later time implies arriving later at vertex $v$. 
The unavailability scenarios observed in TDRPP are addressed by considering unavailability calendars or unavailability intervals. A variant of ARP with unavailability calendars is proposed by \cite{lannez}, where an exact solution is determined using Danzig-Wolfe decomposition and column generation. Lannez's work models this inspection routing and scheduling as an Integer Programming (IP) problem by discretizing time. 
\cite{buriuly2022} models the routing and scheduling problem with unavailability intervals as constraints, using continuous temporal variables. The proposed preliminary column generation approach discusses strategies to compute lower bounds, upper bounds, and sub-optimal solutions for the CARP-TU, and consequently for RPP-TU. The study of exact algorithms is useful for medium-sized networks, as they guarantee optimal solutions, thereby an improvement over sub-optimal solutions, often at the expense of computational runtime. 

\begin{parag}{\textbf{Asked ChatGPT to list research gaps in the previous paragraph} There are no explicit research gaps mentioned in the paragraph. However, some potential research gaps that can be inferred from the paragraph are:}
    \pitemo{The need for a stronger connectivity condition for agent-network compatibility in railway-based routing and scheduling problems with time-varying attributes.}{} 
    \pitemo{The limitation of TDRPP in addressing railway-based routing and scheduling problems due to disconnected sub-graphs caused by unavailability scenarios.}{} 
    \pitemo{The limitation of TDRPP in addressing multi-agent scenarios.}{} 
    \pitemo{The need for more efficient and faster algorithms for solving routing and scheduling problems with unavailability constraints.}{} 
    \pitemo{The potential of continuous temporal variables in modeling routing and scheduling problems with unavailability constraints.}{} 
    \pitemo{The need for exact algorithms to guarantee optimal solutions for medium-sized networks.}{} 
\end{parag} 

 \begin{parag}{What are the research gaps?}
    \pitemo{Intro}{The proposed study on RPP-TRU demonstrates the impact of recurrence in the temporal unavailabilities of RPP-TU. 
    In particular, the novelty of this work is in addressing the incompatibility between temporal requirements of agents and networks; thereby providing criteria for guaranteeing the existence of a solution for RPP-TRU and an exact solution. }
    \pitemo{No study is available to exploit the recurrence property in railway schedules}{Furthermore, unlike the existing works on the planning of periodic train schedules, e.g. \cite{peeters}, no study is observed to exploit the recurrence property of periodic train schedules for routing and scheduling of agents (like inspection wagons) on railway tracks, to the best of our knowledge. }
\end{parag} 

\ifsubbuild

\bibliographystyle{cas-model2-names}

\bibliography{citationsall2.bib}

\fi

\section{Preliminaries }
\label{rpptu:sec:prelims}

\ifsubbuild
\setcounter{section}{2}
\section{Preliminaries}
\fi

Few commonly referred terminologies from graph theory have been revisited in this section. Suitable examples and illustrations are also included throughout the text for clarity. Lastly, a brief overview of the replicated graph is presented for brevity. 

\subsection{Terminologies}
\label{prelims:sec:terms}

Our graph definitions are consistent with \cite{buriuly2022}. Readers may refer to \cite{godsil2001agt} and \cite{akrida2017} for graph theory details, while routing related details for railway problems are available in the survey article by \cite{liden}.

\begin{figure}[!h]
	\centering
     \begin{subfigure}[b]{0.47\textwidth}
         \centering
         \includegraphics[width=\textwidth]{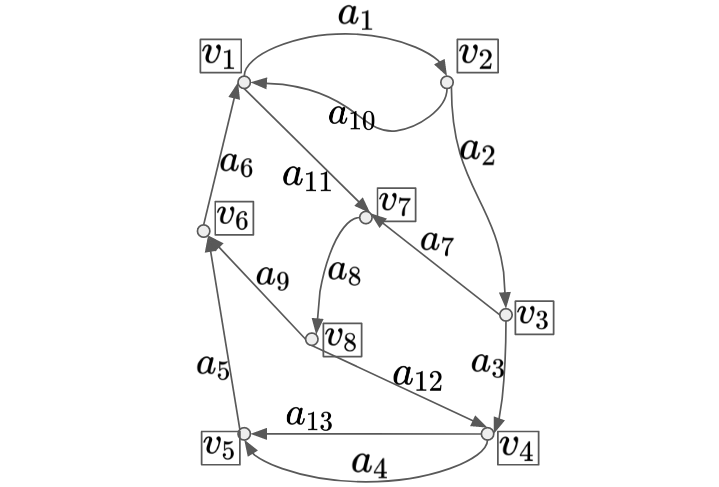}
         \caption{An illustrative graph $G = (V, A, F^+, F^-)$, where $V = \{v_1, \dots, v_8\}$ and $A = \{a_1, \dots, a_{13}\}$, $F^+$ and $F^-$ are suitable maps; source: \cite{buriuly2022}. }
         \label{fig:prelim_fig1a}
     \end{subfigure}
     \hspace{1cm}
     \begin{subfigure}[b]{0.42\textwidth}
         \centering
         \includegraphics[width=\textwidth]{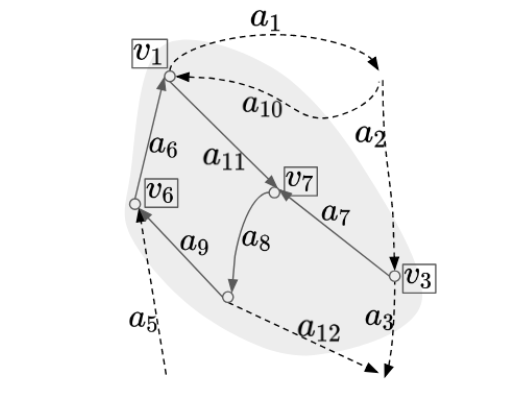}
         \caption{A sub-graph whose vertex set is $\{v_1, v_3, v_6, v_7, v_8\}$ and arc set is $\{a_6, a_7, a_8, a_9, a_{11}\}$. The dashed arcs are \textit{boundary arcs} having either the head or tail vertex within the sub-graph.}
         \label{fig:prelim_fig1c}
     \end{subfigure}
	\captionsetup{justification=centering}
	\caption{A sub-graph is constructed from a multi-graph by selecting all vertices and arcs strictly inside some region (the shaded blob). }
	\label{fig:prelim_fig1}
\end{figure}

\begin{defxy}{Graphs!}{Directed multi-graph}
Directed multi-graphs are useful for representing railway and roadway networks with multiple paths between the same locations. 
Each arc in a directed multi-graph is directed from its tail vertex to its head vertex, and there can be multiple directed arcs between a pair of vertices in a directed multi-graph.
    A directed multi-graph is represented as an ordered 4-tuple $G = (V, A, F^+, F^-)$; where $V$ is a vertex set, $A$ is an arc set, $F^+:A \rightarrow V$ is a map that assigns a head vertex to each arc, and $F^-:A \rightarrow V$ is a map that assigns a tail vertex to each arc.  
    
    Figure \ref{fig:prelim_fig1a} shows an example of directed multi-graph with vertex set $V = \{v_1, \dots, v_8\}$ and arc set $A = \{a_1, \dots, a_{13}\}$, where $a_4$ and $a_{13}$ are parallel arcs. Both arcs $a_4$ and $a_{13}$ are  directed from their tail vertex $F^-(a_4) = F^-(a_{13}) = v_4$ to their head vertex $F^+(a_4) = F^+(a_{13}) = v_5$. 
\end{defxy}

\begin{defy}{Sub-graph}
    A directed graph constructed using a subset of vertices and arcs of graph $G$ is called a sub-graph of $G$. Figure \ref{fig:prelim_fig1c} shows an example of a sub-graph. 
\end{defy}

\begin{defy}{Depot and agent set}
    Depot is the physical start and end location of an \textit{agent} $k \in \mathcal{K}$, represented as one or more vertices of a network graph.
    In general, the vertex $v_1$ is considered as depot for all agents in the \textit{agent set} $\mathcal{K}$. 
\end{defy}

\begin{defy}{Walk, path, and cycle}
    A (directed) \textit{walk} of length $l$ is a sequence of $l$ arcs such that head vertex of every arc (excluding the last arc) in the sequence is a tail vertex of the next arc in the sequence. For example in Figure \ref{fig:prelim_fig1a}, the sequence $a_1, a_{10}, a_{11}, a_8$ is a walk of length $4$, where $F^+(a_1) = F^-(a_{10})$, $F^+(a_{10}) = F^-(a_{11})$, and $F^+(a_{11}) = F^-(a_8)$. If all the vertices in the walk are unique then the walk is called a \textit{vertex-disjoint walk} or \textit{path}. $a_{10}, a_{11}, a_8$ is an example of a path. 

    A \textit{route} is a walk that has the same starting and ending vertex. For example in Figure \ref{fig:prelim_fig1a}, the sequence $a_8, a_9, a_6, a_1, a_{10},a_{11}$ is a route of length $6$, where $F^-(a_8) = F^+(a_{11}) = v_7$. 
    
    A (directed) \textit{cycle} is a vertex-disjoint route. It is also defined as a path that has the same starting and ending vertex. For example in Figure \ref{fig:prelim_fig1a}, the sequence $a_8, a_9, a_6, a_{11}$ is a cycle of length $4$, where $F^-(a_8) = F^+(a_{11}) = v_7$. 
\end{defy} 

\begin{defy}{Service arcs (or required arcs)}
    These arcs are a subset of the arc set of graph $G$, denoted as $A_R \subset A$. A service/required arc represents a physical entity that requires attention, e.g. a rail-track requiring inspection. These arcs must be traversed \textit{at least once}, by \textit{at least one} of the agents. 
\end{defy}

\begin{defy}{Deadhead arcs}
    Arcs that do not require servicing are called deadhead arcs, denoted as $A_D = A\backslash A_R$. The graph composed of all the vertices of graph $G$ and deadhead arcs only is called \textit{deadhead sub-graph}, given as $G_D = (V, A_D, F^+, F^-)$. 
\end{defy}

\begin{defy}{Traversal cost}
    The cost associated with an agent $k \in \mathcal{K}$ for travelling an arc $a_m \in A$ is the traversal cost (or cost of traversal) denoted as $c_{mk}$. Note that each agent $k$ may have different cost while traversing an arc $a_m$.
\end{defy}

\begin{defy}{Running-time}
	Time taken by an agent to traverse an arc is called the running-time (or time of traversal) $t_{mk}$ of the agent $k \in \mathcal{K}$ for the arc $a_m \in A$. For service arcs $a_m \in A_R$, the quantity $t_{mk}$ is termed as \textit{service time}, and for deadhead arcs $a_m \in A_D$, it is termed as \textit{deadhead or cruise time}. The collection of all running-time, also called \textit{running-time data}, is denoted by the matrix $W$, see Figure \ref{fig:prelim_fig2}. 
 
\begin{figure}[!h]
    \centering
    \includegraphics[scale=0.35]{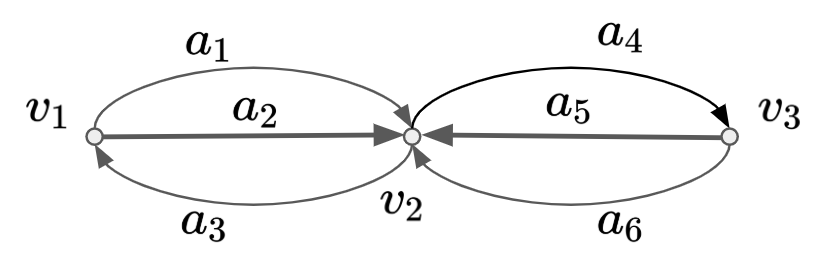}
    \captionsetup{justification=centering}
    \caption{An example graph from \cite{lannez}, where vertex $v_1$ is the depot, and $A_R := \{a_2, a_5\}$ represent the service arcs. The running-time data is given as $W = [1,3,1,1,2,1]^\top$ for single-agent case, and $W = [[1,3,1,1,2,1]^\top,[1,3,1,1,2,1]^\top]$ for two-agent case. }
    \label{fig:prelim_fig2}
\end{figure}
\end{defy}

\begin{defy}{Strongly connected graph}
    A  directed graph is said to be strongly connected if there exists a pair of paths joining any two vertices in either direction.
\end{defy}

\begin{defy}{Unavailability list}
    The unavailability list (or the list of temporal unavailabilities) for arc set $A$ is denoted as $Z = \cup_{a_m \in A} Z_m$. Each entry $(z^-, z^+) \in Z_m$ is a two-tuple that describes the time interval of duration `$z^+-z^-$' between which arc $a_m$ is unavailable. 
\end{defy} 

\begin{defy}{Temporal graph-agent tuple}
    We define a temporal graph-agent tuple as $3$-tuple $(G,Z,W)$; where $G$ is a multi-graph representing some network, $Z$ is the list of unavailabilities for arc set $A$ of graph $G$, and the vector $W$ is the collection of running-times for all agents. 
\end{defy} 

\begin{defy}{$T_p$-recurrence}
    A temporal graph-agent tuple $(G,Z,W)$ is $T_p$-recurrent if each arc in the parent graph $G$ is available at least once in $T_p$-period interval. In particular, a snapshot of the graph might appear different for different time instants as some arcs are unavailable at certain time intervals, however, $T_p$-recurrence ensures that overlaying all these graph instances for any arbitrary period $T_p$ will result in the parent graph $G$. 
    The definition is synonymous to $\delta$-recurrence for temporal graph, see \cite{ilcinkas2018}. 
\end{defy}

\begin{defy}{Cutting-planes or cuts} 
    Cutting-planes or cuts are inequality constraints that refine the feasible set e.g. removing non-integer solutions by adding an inequality constraint to the formulation of an integer programming problem. 
\end{defy}

\subsection{Replicated graph}
\label{prelims:sec:repg}

The concept of a replicated graph is 
a graphical extension to the 
Arc Path Arc Sequence (APAS) formulation presented in \cite{tan11a}, which models the order of traversal using layered representation. The replicated graph is useful in representing temporal/time-varying data. We present a brief description of the replicated graph introduced in \cite{buriuly2022} for completeness.

\begin{figure*}[!h]
	\centering
	\includegraphics[width=0.61\textwidth]{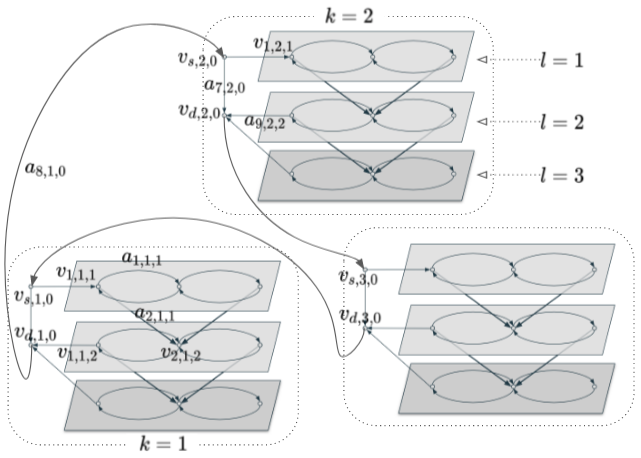}
	\captionsetup{justification=centering}
	\caption{A replicated graph for $|\mathcal{K}|$ agents, constructed from a base graph shown in Figure \ref{fig:prelim_fig2}. }
	\label{fig:prelim_fig3}
\end{figure*}

The replicated graph $\mathcal{G} = (\mathcal{V}, \mathcal{A}, \mathcal{F}^+, \mathcal{F}^-)$ is a multi-graph derived from the base graph $G = (V, A, F^+, F^-)$. The running-time data for each arc of the replicated graph is denoted using vector $\mathcal{W}$. Figure \ref{fig:prelim_fig3} shows a replicated graph with $|\mathcal{L}|$ layers\footnote{It is a copy of the deadhead sub-graph $G_D$. } (illustrated using planes) in each of the $|\mathcal{K}|$ \textit{agent-sub-graph} 
(illustrated using dotted boxes). In this figure, $k \in \mathcal{K} (= \{1,2,3\})$ is the agent set, that dictates the number of agent-sub-graphs. Given required arc set $A_R = \{a_2, a_5\}$, each agent-sub-graph has $|A_R|+1$ layers, denoted as $l \in \mathcal{L} (= \{1,2,3\})$. The deadhead graph $G_D$ is obtained from $G$ in Figure \ref{fig:prelim_fig2}, and replicated in each layer of $\mathcal{G}$. The required arcs in $A_R$ are used to connect the layers in increasing order of the layer numbers/identifiers. Additionally, source and destination vertices (indexed as `$s$' and `$d$' respectively) are introduced for each agent $\mathcal{K}$ and connected to the copies of the depot vertex in each layer. The source vertex connects to the depot of the first layer, while the destination vertex is connected to the depot vertices of the rest of the layers, to ensure the correct flow of decisions. 

In the replicated graph, the vertices are represented with three index-subscripts separated by commas, e.g. $v_{i,k,l} \in \mathcal{V}$, where $k \in \mathcal{K}$, $l \in \{0\} \cup \mathcal{L}$, and $v_i$ is either a vertex of base graph $G$ (if $l \geq 1$) or an extra vertex (with index `$s$' and `$d$'; and layer $l = 0$). For brevity, \textit{all similar notations with index-subscripts are often written without a comma separation}, e.g. $v_{ikl}$. Similar convention is followed for denoting the arcs of replicated graph; e.g. arc $a_{2,1,1} \in \mathcal{A}$ in Figure \ref{fig:prelim_fig3} is a copy of arc $a_2 \in A$ from Figure \ref{fig:prelim_fig2} placed on layer $l=1$ and agent-sub-graph $k=1$. Also, note that the agent-sub-graphs are interconnected to form a cycle connecting source and sink vertices in alternate: $\{v_{s,2,0},a_{7,2,0},v_{d,2,0},\dots,v_{d,3,0},\dots,v_{d,1,0},a_{8,1,0},v_{s,2,0}\}$\footnote{Ideally a cycle in multi-graph is described using arcs only $\{a_{7,2,0},\dots,a_{8,1,0}\}$, however, here we have also included the vertex notations in alternate for brevity.}. This cycle is comprised of three \textit{source-destination arcs} and three \textit{intra-layer arcs} in the $0^{th}$ layer. In Figure \ref{fig:prelim_fig3}, the arc $a_{7,2,0}$ is a source-destination arc and $a_{8,1,0}$ is an intra-layer arc, where the arcs $a_7$ and $a_8$ are virtual arcs i.e. not present in the base graph $G$. 
    
In each layer of the replicated graph, copies of the base graph vertex set $V$ and the deadhead arc set $A_D$ are embedded, resulting in $|{V}||\mathcal{K}||\mathcal{L}|$ vertices and $|{A_D}||\mathcal{K}||\mathcal{L}|$ arcs. 
There are $2|\mathcal{K}|$ additional vertices to represent source and destination (vertices indexed with $s$ and $d$), and $(|\mathcal{K}||\mathcal{L}|+2)$ extra arcs adjacent to these vertices. 
$|\mathcal{K}|(|\mathcal{L}|-1)$ copies of the service arcs $A_R$ of the base graph $G$ are introduced between the layers such that an agent can only traverse in increasing order of the layers (top-to-bottom). 
Due to this layered construction, traversing two different arcs of replicated graph implies re-traversing the same arc of the base graph. However, in a replicated graph, the order of traversing this arc is unambiguous because one with a smaller layer index is traversed first, reflected by the top-to-bottom arc directions. 

The \textit{cost of traversal} and \textit{time of traversal} for arcs $a_{mkl} \in \mathcal{A}$ of replicated graph $\mathcal{G}$ are denoted as $c_{mkl}$ and $t_{mkl}$ respectively. Note that, $c_{mkl} = c_{mk}$ and $t_{mkl} = t_{mk}$, for all $l \in \mathcal{L}$. 

\subsection{Arc unavailability to vertex inactivity} 
\label{prelims:sec:arc2vxunav} 

\begin{figure}[!h]
    \centering
	\includegraphics[width=0.65\textwidth]{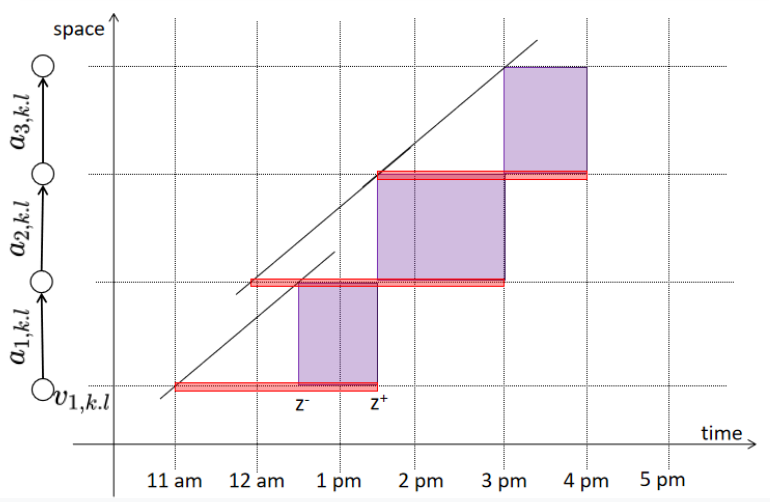}
	\captionsetup{justification=centering}
	\caption{Three arcs of a graph $G$ having different unavailability periods are shown. Based on the running-time data $W$ (illustrated with slanting lines) and the unavailability list $Z$ (fat blocks), the infeasible regions for the departure time of an agent (narrow boxes over dashed lines) at each of the tail vertices of corresponding arcs are computed. This restriction guarantees that for any arc of the replicated graph, if an agent's departure time at the tail vertex is outside the highlighted infeasible region, then the agent will reach the other end without entering the unavailability zone (fat blocks) of the arc and then wait as long as required at the tail vertex of next arc. For instance, if an agent departs at 10:59 am from the tail vertex of arc $a_1$, denoted as $v_{1,k,l}$ then the agent arrives at the head vertex before the arc gets blocked at 12:29 pm. The agent then waits till 3:00 pm as the next arc $a_2$ will be blocked midway of its traversal, if it departs at 12:29 pm. }
	\label{fig:unav}
\end{figure}

Restricting the agents to access an unavailable arc is equivalently achieved/modeled by prohibiting agents' activities at the vertices under the assumption that agents are not allowed to stop midway on an arc. In particular, implementing the waiting of agents only at the tail vertex of an arc is sufficient to model the unavailability restrictions. Figure \ref{fig:unav} illustrates the departure time restrictions at the vertex due to arc unavailabilities. Given a replicated graph $\mathcal{G}$ and a list of temporal unavailabilities $Z$, the time value at vertex $v_{1,k,l} \in \mathcal{V}$ should be either less than $z^- - t_{1,k,l}$ or more than $z^+$; where $v_{1,k,l}$ is the tail vertex of arc $a_{1,k,l}$, $(z^-,z^+) \in Z_1$ indicate an unavailability period of arc $a_{1,m,l}$, and the time of traversal for arc $a_{1,k,l}$ is $t_{1,k,l} = 1\textit{hr } 30\textit{mins}$ in Figure \ref{fig:unav}. Multiple unavailability restrictions lead to disjoint feasible sets, thus a set of disjunctive linear constraints over the departure times of agents at the vertices is sufficient to model waiting decisions (see Equation \eqref{form:unav} in Section \ref{rpptu:sec:formulate}). Without loss of generality, waiting for the availability of an arc is executed at its tail vertex only. 

\ifsubbuild

\bibliographystyle{cas-model2-names}

\bibliography{citationsall2.bib}

\fi

\section{Rural Postman Problem with Temporal Recurring Unavailabilities}
\label{rpptu:sec:formulate}

\ifsubbuild
\setcounter{section}{2}
\section{RPP-TRU}
\fi

\begin{parag}{Whats here?}
 \pitemo{About}{In this section, we discuss the formal problem statement for RPP-TRU along with a Mixed Integer Linear Programming (MILP) formulation for determining routing and scheduling decisions. A three-index formulation for RPP-TRU is constructed using a replicated graph, obtained from the directed multi-graph $G$ in the RPP-TRU problem statement. 
 We also study the effect of recurrence in RPP-TRU, and analyze the underlying network.} 
 %
\end{parag}

\subsection{Problem statement }


RPP-TRU is described over a \rindxy{Directed multi-graph}{directed multi-graph} $G = (V,A,F^+,F^-)$ that represents some network, e.g. a portion of railway network. The objective of RPP-TRU is to determine routes and schedules for $|\mathcal{K}|$ agents such that the total cost of traversal and time of task completion (maximum among the time of traversal of all agents) is minimized. The trade-off between the spatial (traversal) and temporal (completion time) costs is dependent on the practical problem under study. The following constraints must be satisfied by the solution plans: the routes must collectively service a subset of arcs called service/required arcs $A_R \subset A$, while the schedule must respect the running-time data $W$ and the list of unavailabilities $Z$. The unavailabilities recur at different intervals, e.g in railways, the unavailabilities are based on the schedule of trains and track maintenance possessions\footnote{\textit{Track possession} means reserving a track section for some activity (inspection, maintenance, etc) for a given period of time \cite{liden}. }. 

\subsection{Formulation }
\label{rpptu:sec:theformulation}

\keyq{How do we formulate the problem, and why is it valid?} 
We consider the \textit{temporal graph-agent tuple} $(G,Z,W)$ as \textit{$T_p$-recurrent}, where $G$, $Z$ and $W$ denote the network, unavailability and running-time respectively. 
In this section, we present a three-index formulation of the RPP-TRU (similar to RPP-TU from \cite{buriuly2022}) based on the replicated graph $\mathcal{G}$ (constructed from $G$), running-time of a replicated graph $\mathcal{W}$, and the unavailability list $Z$. The RPP-TRU is formulated with a generalized cost term to penalize both the spatial movements and the time of completion. 
Adapting the three indices from the \textit{replicated graph} notations ensures that the routing portion of the proposed mathematical formulation is directly derived by applying the commonly known flow-conservation based approach on the replicated graph. The scheduling portion of the problem is also benefited from the three-index formulation as the order of traversal becomes intrinsic. The proposed formulation also ensures that the assignment of only one temporal variable at each vertex of the replicated graph is sufficient to model the optimization problem (see Corollary \ref{cor:sufficient} in Section \ref{sec:prop}), provided it holds some connectivity properties. 

\subsubsection{Decision variables } 
\label{rpptru:sec:decvar}

Let $X \in \{0,1\}^{|\mathcal{A}|}$ and $\Gamma \in \mathbb{R}^{|\mathcal{V}|}, \Gamma \geq 0$ be the spatial and temporal decision vectors respectively. In indexed representation\footnote{Instead of representing the three-indexed spatial and temporal decision variables as $3$D array we represent them as $1$D vector. This allows for matrix multiplication with the incidence matrix in later sections. } , the spatial routes of agents are indicated by $X_{mkl}$ for every arc $a_{mkl} \in \mathcal{A}$, while the schedule decisions are indicated by $\Gamma_{ikl}$ for every vertex $v_{ikl} \in \mathcal{V}$ of the replicated graph 
$\mathcal{G}$. Here, $m$ and $i$ are identifiers for arc $a_m \in A$ and vertex $v_i \in V$ respectively, $k \in \mathcal{K}$ is a positive integer indicating an agent, and $l \in \mathcal{L}$ is a non-negative integer indicating the order of traversal. Other relevant identifiers are $m', j, k', l'$, where $a_{m'} \in A$, $v_j \in V$, $k' \in \mathcal{K}$ and $l' \in \mathcal{L}$. 
\begin{equation*}
 \nonumber 
 \begin{aligned}
  X_{mkl} &= && \begin{cases} 1, & \text{if arc \hgl{$a_{mkl}$ is traversed/served i.e. arc }$a_m$ is traversed/served by agent $k,$ after servicing the }\\
  & \text{$(l-1)^{th}$ required arc and before servicing the $l^{th}$ required arc}  \\ 
  0, & \text{otherwise} \end{cases} \\ 
  \Gamma_{ikl} &\rightarrow && (\geq 0), \text{ departure time\hgl{ at vertex $v_{ikl}$, or departure time} of agent $k$ at vertex $v_i$, after servicing the }\\ 
  & && \text{$(l-1)^{th}$ required arc and before servicing the $l^{th}$ required arc} 
 \end{aligned}
\end{equation*}

\subsubsection{Three-index formulation} 
\label{rpptru:sec:3inf}


The three-index formulation for RPP-TRU, denoted as \ref{form:rpptu3if}, is composed of: (i) spatio-temporal cost in Equation \eqref{form:rpptucost}, (ii) spatial constraints (denoted as $\mathcal{S}_X$ in brief) described using Equations \eqref{form:flow}-\eqref{form:service}, and (iii) spatio-temporal constraints described using Equations \eqref{form:maxtime}-\eqref{form:unav}. 
The trade-off between the spatial cost $\sum_{m, k, l} c_{mkl} X_{mkl}$ and temporal cost $\gamma$ is determined by the factor $\beta$. Table \ref{tab:formnotas} summarizes all important notations used in the formulation. 

\begin{table}[!htp]\centering
\caption{Review of important notations used in this section. }
\label{tab:formnotas}
\begin{tabular}{c|p{10cm}|c}\toprule
Notations & Description & Refs.\\ 
\hline 
$G$ & $(:= (V, A, F^+, F^-))$ \rindxy{Directed multi-graph}{Multi-graph} from the problem statement, that represents or models a practical network. & Defined in Sec-\ref{prelims:sec:terms} \\ 
$G_D$ & $(:= (V, A \backslash A_R, F^+, F^-))$ Multi-graph with service/required arcs removed, known as deadhead graph or \rindxy{deadhead sub-graph}{deadhead sub-graph}. & Defined in Sec-\ref{prelims:sec:terms} \\ 
$\mathcal{G}$ & $(:= (\mathcal{V},\mathcal{A},\mathcal{F}^+,\mathcal{F}^-))$ \rindx{Replicated graph} obtained from the multi-graph $G$ and agent set $\mathcal{K}$. & Defined in Sec-\ref{prelims:sec:repg} \\ 
$\mathcal{L}$ & Layer set: to indicate the order of servicing as an intrinsic part of the replicated graph. The elements of the layer set are indicated by either $l$ or $l'$. & Defined in Sec-\ref{prelims:sec:repg} \\ 
$W, \mathcal{W}$ & Running-time data in matrix form: $W$ is composed of running-time entries $t_{mk}$, and $\mathcal{W}$ is composed of running-time entries of replicated graph $t_{mkl}$. & Defined in Sec-\ref{prelims:sec:terms} and \ref{prelims:sec:repg} \\ 
$Z, Z_m$ & The list of unavailability intervals of arc $a_m \in A$ of graph $G$ is denoted by $Z_m$. Each entry in $Z_m$ is comprised of a lower time limit $z^-$ and an upper time limit $z^+$ to denote the unavailability interval. All the unavailability lists are collectively represented by $Z$. & Defined in Sec-\ref{prelims:sec:terms} \\ 
$\delta^-_{ikl}, \delta^+_{ikl}$ & Denotes the set of incoming and outgoing arcs to and from the vertex $v_{ikl} \in \mathcal{V}$ of the replicated graph $\mathcal{G}$: Mathematically, $\delta^-_{ikl} = \{a_{mkl'} | \mathcal{F}^+(a_{mkl'}) = v_{ikl}\}$ and $\delta^+_{ikl} = \{a_{mkl'} | \mathcal{F}^-(a_{mkl'}) = v_{ikl}\}$; where $l$ maybe different from $l'$. & Used in Eq-\eqref{form:flow} \\ 
$\mathcal{A}_m$ &  All copies of a service arc $a_m \in A_R$ in the replicated graph $\mathcal{G}$ is given as: $\mathcal{A}_m = \{a_{mkl} \in \mathcal{A}_R \ | ~\mbox{where }k \in \mathcal{K}, \mbox{ and }l \in \mathcal{L}\}$. We also have the set of required arcs as a disjoint union of these copies: $ \cup_{a_m \in A_R} \mathcal{A}_m = \mathcal{A}_R$. & Used in Eq-\eqref{form:service} \\ 
$X, \Gamma$ & Denotes the spatial and temporal decision variables (vectors) respectively. Corresponding index notations are $X_{mkl}$ and $\Gamma_{ikl}$/$\Gamma_{jkl}$.  & Defined in Sec-\ref{rpptru:sec:decvar} \\ 
$\beta$ & Trade-off parameter between penalty on movement and delay. & Defined in Sec-\ref{rpptru:sec:3inf} \\ 
$c_{mkl}, t_{mkl}$ & Cost and time of traversal of arc $a_{mkl} \in \mathcal{A}$. & Defined in Sec-\ref{prelims:sec:repg} \\ 
$\mathcal{S}_X$ & Denotes the spatial constraints in the formulation of RPP-TRU: described using Equations \eqref{form:flow}-\eqref{form:service}. These constraints depend only on the spatial variables $X$. & - \\ 
\bottomrule
\end{tabular}
\end{table}

\setcounter{eqc}{\theequation}
\stepcounter{eqc}
\begin{subequations}
\makeatletter
\def\@currentlabel{P\arabic{eqc}}
\makeatother
\label{form:rpptu3if}
\renewcommand{\theequation}{P\arabic{eqc}-\alph{equation}}
 \begin{alignat}{1}
	\label{form:rpptucost}
 	\min_{X, \Gamma, \gamma} ~&\sum_{m, k, l} c_{mkl} X_{mkl} + \beta ~ \gamma \\ 
	\label{form:flow}
 	\textit{s.t.} ~& 
    \sum_{a_{mkl} \in \delta^+_{ikl}} X_{mkl} - \sum_{a_{mkl}  \in \delta^-_{ikl}} X_{mkl} = 0, ~\forall v_{ikl} \in \mathcal{V} \\
	\label{form:service} 
	& \sum_{a_{mkl} \in \mathcal{A}_m} X_{mkl} = 1,\ \ \forall a_m \in A_R   \\
	\label{form:maxtime}
 	& \gamma \geq \Gamma_{dk0}, ~ \forall k \in \mathcal{K} \\ 
	\label{form:run} 
	& \Gamma_{jkl'} \geq t_{mkl} + \Gamma_{ikl} + \tau(1 - X_{mkl}) ; ~\forall a_{mkl} \in \mathcal{A} \\ 
    \nonumber 
    &~~~~ \Gamma_{ikl} \textit{ is the temporal value at tail vertex }v_{ikl} = \mathcal{F}^-(a_{mkl}), \\ 
    \nonumber 
    &~~~~ \Gamma_{jkl'} \textit{ is the temporal value at head vertex }v_{jkl'}=\mathcal{F}^+(a_{mkl}) \\ 
	\label{form:unav}
 	& \left.
    \begin{array}{ll}
     either,&\\
     ~~\Gamma_{jkl} \leq &z^- - t_{mkl} + \tau (1 - X_{mkl})\\
     or, &\\
     ~~\Gamma_{jkl} \geq &z^+ - \tau (1 - X_{mkl}) \\
    \end{array}
~~~~~~ \right\} \forall (z^-, z^+) \in Z_m, a_{mkl} \in \mathcal{A} \\ 
    \nonumber 
     &~~~~ \textit{where,} \Gamma_{jkl} = F^-(a_{mkl}) \\ 
	\label{form:intlin}
 	& X_{mkl} \in \{0,1\}, \Gamma_{jkl} \geq 0 
 \end{alignat}
\end{subequations} 

The expression \eqref{form:rpptucost} imposes minimum traversal cost and minimum time of completion of all tasks. The factor $\beta$ determines the trade-off between the spatial and temporal cost terms. Each routing decision in $X_{mkl}$ is penalized by traversal cost $c_{mkl}$, while the term $\gamma$ (by design) attains the maximum among completion time of all agents due to Equation \eqref{form:maxtime} i.e. smallest value of $\gamma$ equals $\max_{k\in \mathcal{K}} \Gamma_{dk0}$. Equation \eqref{form:flow} ensures that the spatial solutions are valid by treating decisions as flow in the replicated graph. These flow constraints balance the incoming and outgoing decisions at vertices. Equation \eqref{form:service} assigns each service arc to only one of the agents such that only one among all copies of a service arc $\mathcal{A}_m$ must be serviced. \hgl{The running-time constraints}\footnote{They also eliminate sub-tours from the spatial decisions, see \cite{mtz}.} \hgl{in Equation \eqref{form:run} determines a schedule that respects the travel time of the arcs, thereby managing the time difference between the head and tail vertex of an arc $a_{mkl} \in \mathcal{A}$. These constraints are imposed only when the spatial decision $X_{mkl}$ is $1$ for arc the $a_{mkl}$ i.e. if $X_{mkl}=0$ then the right-hand side of the constraint is a very large negative number $-\tau (1 - X_{mkl}) = -\tau$. If $X_{mkl}=1$ then the time difference at the head and tail vertex $\Gamma_{jkl'} - \Gamma_{ikl}$ must be larger than a positive number $t_{mkl}$. Lastly, the unavailability constraints in Equation \eqref{form:unav} compels the agents to wait at the vertices for the availability of the outgoing arcs. In particular, if $X_{mkl} = 1$ then the unavailability constraints ensure that $\Gamma_{jkl}$ is not between $z^- - t_{mkl}$ and $z^+$, otherwise $\Gamma_{jkl}$ is upper bounded by a large number $O(\tau)$.} The resulting time value $\Gamma_{jkl}$ is either smaller than the lower limit $z^--t_{mkl}$ or larger than the upper limit $z^+$ provided the outgoing arc is part of the spatial solution; see Section \ref{prelims:sec:arc2vxunav}. Note $\tau$ must be very large i.e. an upper bound on $\Gamma$, see Corollary \ref{cor:timebound} for choosing suitable $\tau$. In the proposed algorithm, $\tau$ is eliminated automatically, hence we don't need an explicit value for it. $\tau$ is useful only when directly solving the RPP-TRU using MILP solvers. 

In order to apply a MILP solver, the disjunctive constraints in Equation \eqref{form:unav} is augmented with dummy binary variables $\chi_{\sigma}$ as shown in Equation \eqref{eqn:unav_dummy}. Since $\chi_{\sigma} \in \{0,1\}$, only one of the two constraints will be active. In the proposed algorithm, the unavailability constraints will be implemented using a graph search algorithm, to avoid introducing a large number of dummy variables. 

\begin{equation}
 \label{eqn:unav_dummy} 
 \begin{aligned}
 	& \left.
    \begin{array}{ll}
     ~~ \Gamma_{jkl} \leq &z^- - t_{mkl} + \tau (1 - X_{mkl}) + \tau \chi_{\sigma}\\
     ~~ \Gamma_{jkl} \geq &z^+ - \tau (1 - X_{mkl})  - \tau (1-\chi_{\sigma})\\
     \textit{where,} & \textit{the } \sigma^{th} \textit{ entry } (z^-, z^+) \in Z_m \\ 
     & \textit{and } \Gamma_{jkl} = F^-(a_{mkl}), \textit{ tail vertex or the arc}\\ 
    \end{array}
\right\} ~~\forall \sigma = 1, \dots, |Z_m|; a_{mkl} \in \mathcal{A} 
 \end{aligned}
\end{equation}

\subsection{Properties of the replicated graph based formulation}
\label{sec:prop}

We first ensure that the RPP-TRU problem dataset is \textit{well-defined} in the upcoming discussions. Next, we establish various results for the replicated graph based RPP-TRU formulation, based on well-defined dataset and $T_p$-recurrence. 

\noindent A temporal graph-agent tuple $(G,Z,W)$ (see definition in Section \ref{prelims:sec:terms}) is said to be \textit{well-defined} if \hgl{all of the following conditions are satisfied}:
\begin{itemize}
	\item[1.] every arc $a_m \in A$ in $G$ is available for at least the corresponding running-time $t_m (:= \min_{k \in \mathcal{K}} t_{mk})$ for any one agent $k \in \mathcal{K}$, in the period $[0,\infty)$, 
	\item[2.] availability periods\footnote{time gap between two unavailability periods}, derived from $Z$, for every arc in $G$ \hgl{that are smaller than the running-time of that arc for all the agents, are removed from $Z$ }
    \hgl{(small availability gaps are equivalent to being} unavailable for the agents), and
	\item[3.] the deadhead graph $G_D := (V, A_D, F^+, F^-)$ (see Table \ref{tab:formnotas}), is strongly connected i.e. there exists a path between any two vertices in the deadhead graph. 
\end{itemize} 

\hgl{The spatial connectivity of the network is guaranteed by condition-3 for the well-defined temporal graph-agent tuple. In order to satisfy condition-2, all availability gaps that are completely non-traversable by all the agents are removed from $Z$. This may, however, result in arcs that are never available. Therefore, condition-1 ensures that an arc is available at least once for at least one agent in the period $[0, \infty)$. Consequently, every arc is guaranteed to be traversable if its tail vertex is the depot i.e. an agent takes zero time to reach the arc. For instance, an arc with running-time of one unit might only be available for interval $[0,1]$, and hence it is not traversable by any agent that is not present at the tail vertex of this arc at $t=0$. In conclusion, a stronger temporal condition is still required for guaranteed temporal connectivity of the network. }

Note that, it is possible to construct an example for RPP-TU/RPP-TRU without any solution if the temporal graph-agent tuple is ill-defined. \hgl{For instance, very small availability gaps in the unavailability list $Z$ will inhibit agents from traversing these arcs. This example is indicative of bad agent-network pairing i.e. deploying slow moving agents in busy networks.} Here, we make a stronger statement in Theorem \ref{thm:wd} that guarantees the existence of an optimal solution in the RPP-TU/RPP-TRU formulation, under reasonable assumptions. 

\begin{thm}
\label{thm:wd} 
The replicated graph $\mathcal{G}$ with $|A_R|+1$ layers, contains an optimal solution for the corresponding RPP-TRU in \ref{form:rpptu3if}, if the underlying temporal graph-agent tuple $(G, Z, W)$ is well-defined and $T_p$-recurrent ($T_p \geq 0$), and the cost of traversal $c_{mk}$ is non-negative $\forall a_k \in A, ~k \in \mathcal{K}$. 
\end{thm}

\begin{prf}
Every layer of the replicated graph $\mathcal{G}$ is constructed from deadhead graph $G_D$, hence we can map walks in $G_D$ to an equivalent walk in $\mathcal{G}$. For simplicity in notations, a sequence is constructed in $G_D$ as $\{v_{a_1},a_1,v_{a_2}, a_2, \dots, a_r, v_{a_{r+1}}\}$, with vertices \authar[$v_{a_1},v_{a_2},\dots,v_{a_{r+1}} \in V$ and arcs $a_1,a_2,\dots,a_r \in A\backslash A_R$]. Let $a_1, a_{r+1} \in A_R$ be two arbitrary service arcs; hence this sequence is assumed to be a walk between the tail vertex $v_{a_1}$ of arc $a_1$ and head vertex of arc $a_{r+1}$ (denoted as $v_{a_{r}}$). 
Since every arc is traversable in any $T_p$ time interval - an implication of being well-defined $(G,Z,W)$ and $T_p$-recurrent - feasible time-stamps can be assigned to all vertices of the walk making it a feasible trajectory with respect to unavailabilities. If this walk forms a cycle within a layer, then all arcs of the cycle can be contracted to form a single vertex whose time-stamp is the maximum among all time-stamps in the cycle, thus preserving feasibility by ensuring incremental time stamps. Hence, \authar[the walk $P (:=\{v_{a_1},a_1,v_{a_2}, a_2, \dots, a_{r-1}, v_{a_{r}}\})$] can be reduced to a path (vertex-disjoint walk) preferring waiting instead of deadhead cyclic traversal. 

Given $c_{mkl} \geq 0$, the cost of per unit waiting time is at most the cost of per unit deadhead traversal time. Hence, there exists \authar[such] an optimal path between any two vertices in the deadhead graph $G_D$ i.e. in every layer of the replicated graph $\mathcal{G}$. Therefore given the order of attending the service arcs in an optimal solution, one can construct an optimal trajectory by joining the source vertices (marked with subscript $s$), the destination vertices (marked with subscript $d$), and the service arcs, using suitable paths. In particular, an optimal solution of RPP-TRU is present in our well-defined and $T_p$-recurring formulation using the replicated graph, with dataset governed by graph $G$, unavailability $Z$, running-time data $W$, and traversal cost data $c_{mkl}, ~\forall a_m \in A, k \in \mathcal{K}$.
\end{prf} 

\hgl{Note that, $T_p$-recurrent condition in Theorem }\ref{thm:wd}\hgl{ ensures that condition-1 of well-defined $(G,Z,W)$ is true at every $T_p$ interval, i.e. every arc $a_m$ is traversable in the interval $[t,t+T_p+t_{mk}]$ for some agent $k \in \mathcal{K}$, where $t$ is any arbitrary time in $[0, \infty)$. Next, }Corollary \ref{cor:sufficient} and \ref{cor:timebound} are obtained from Theorem \ref{thm:wd}. 
These results aid in the selection of $\tau$ in Equations \eqref{form:run}, \eqref{form:unav}, and \eqref{eqn:unav_dummy}. 

\begin{corr}
\label{cor:sufficient}
    The optimal solution of RPP-TRU is a cycle in the replicated graph $\mathcal{G}$ that doesn't visit any vertex more than once\hgl{ (i.e. one incoming and one outgoing arc per vertex)}, implying that a single temporal variable per vertex is sufficient to describe the schedules. 
\end{corr}

\begin{corr}
\label{cor:timebound}
    An upper bound on temporal value $\Gamma_{i,k,l}$ is given as $(|V||\mathcal{L}|+1)(T_p + \max_{a_{mkl} \in \mathcal{A}} t_{mkl})$. 
\end{corr}

\begin{prf}
 In each agent-sub-graph of the replicated graph, there are $|V||\mathcal{L}| + 2$ vertices (see Section \ref{prelims:sec:repg}), therefore any path will have a maximum length of $|V||\mathcal{L}| + 1$. Each path is available \hgl{for traversal by at least one of the agents in the period $T_p$ due to periodic connectivity, and the maximum waiting + traversal time is} $T_p + \max_{a_m \in A, k \in \mathcal{K}} t_{m,k}$. 
\end{prf}

Theorem \ref{thm:wd} also establishes Remark \ref{rem:feasball}, which elucidates that the unavailability constraints don't effect the spatial solution \hgl{$X \in \mathcal{S}_X$ (see Table }\ref{tab:formnotas}\hgl{). The reason being every temporal solution can be modified based on the unavailability of arcs, by waiting at the corresponding tail vertex of the arc till the arc is available. All arcs are guaranteed to be available in the near future due to well-defined and $T_p$-recurrence attributes. }

\begin{rem}
\label{rem:feasball} 
Removal of unavailability constraints doesn't alter the spatial feasible region given by Equations \eqref{form:flow} and \eqref{form:service}. 
\end{rem} 

\ifsubbuild

\bibliographystyle{cas-model2-names}

\bibliography{citationsall2.bib}

\fi

\section{Algorithm and analysis}
\label{rpptu:sec:algo}

\ifsubbuild
\setcounter{section}{4}
\section{Algorithm}
\fi

In this section, we propose two approaches for implementing the modified Benders' algorithm for RPP-TRU. Next, we discuss the modification process to arrive at the proposed Benders' decomposition algorithm. Relevant theoretical analysis is also presented in this section to inherit the convergence property of the well-established Benders' decomposition algorithm. 

\subsection{Proposed algorithms}
\label{rpptu:sec:propbendalgo}

\begin{figure}[!h]
	\centering
     \begin{subfigure}[b]{0.89\textwidth}
         \centering
    	\includegraphics[width=0.98\textwidth]{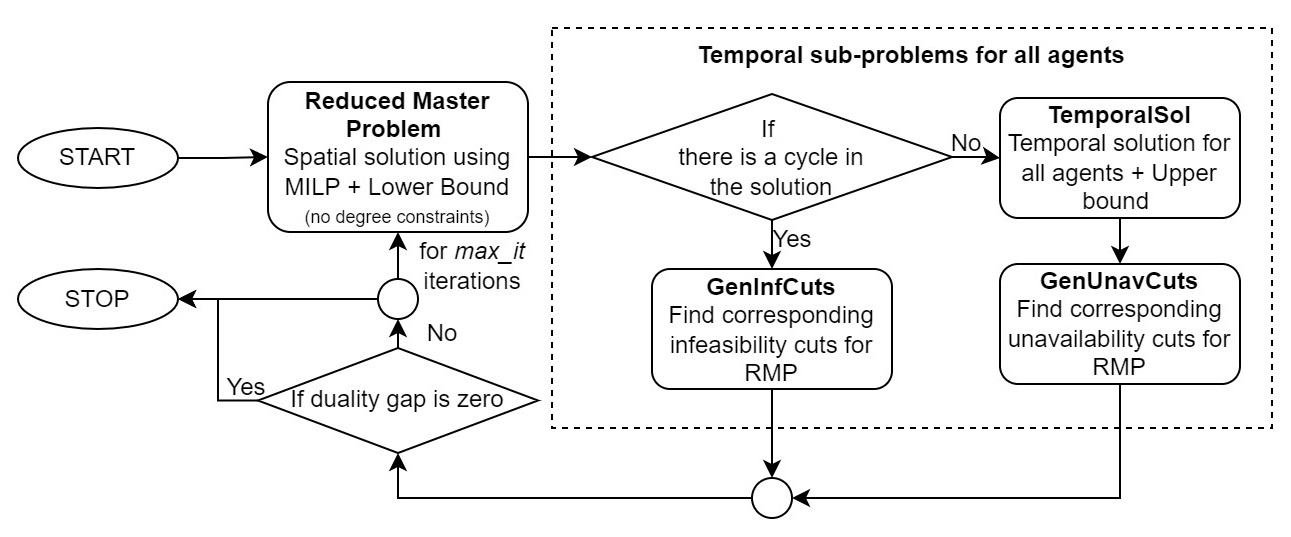}
         \caption{First approach. }
         \label{fig:bndrs1_rpp-tu}
     \end{subfigure}
      \hspace{1cm}
     \begin{subfigure}[b]{0.64\textwidth}
         \centering
    	\includegraphics[width=0.98\textwidth]{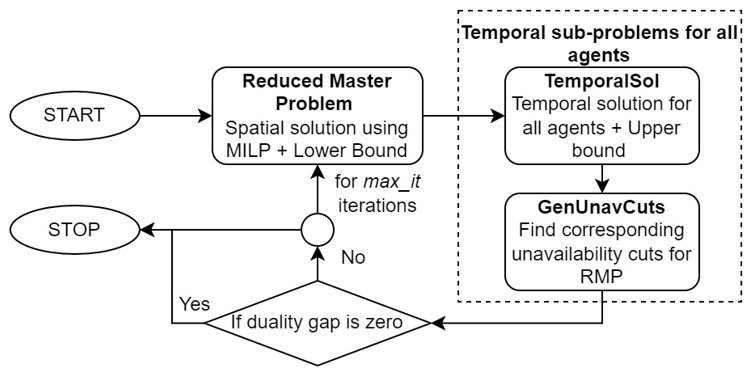}
         \caption{Second approach.}
         \label{fig:bndrs2_rpp-tu}
     \end{subfigure}
	\captionsetup{justification=centering}
	\caption{Proposed Benders' decomposition algorithm variations for RPP-TRU, showing reduced master problems \ref{form:rmp1} and \ref{form:rmp2}, and temporal sub-problems \ref{form:bsp} for each $k \in \mathcal{K}$. The pseudo-code for \textbf{TemporalSol} and \textbf{GenUnavCuts} is shown in Algorithm \ref{algo:temporal} and \ref{algo:unav} respectively.  }
	\label{fig:bndrs_rpp-tu} 
\end{figure}

In Figure \ref{fig:bndrs_rpp-tu} we describe the flowchart of two modified variations of the Benders' decomposition algorithm. Both the algorithms are composed of two main steps: the reduced master problem, and the temporal sub-problems. Furthermore, the temporal sub-problems are divided into two/three computation blocks: \textbf{TemporalSol}, \textbf{GenInfCuts}, and \textbf{GenUnavCuts}, as shown in Figure \ref{fig:bndrs_rpp-tu}. The reduced master problem generates a spatial solution and a lower bound to the optimal solution, while the temporal sub-problems collectively generate an optimal temporal solution corresponding to the spatial solution. When the spatial and temporal solution pair is feasible for RPP-TRU, the corresponding cost value gives an upper bound for RPP-TRU. 
In case there exists an intra-layer cycle in the spatial solution, there is no corresponding temporal solution for the given spatial solution. Hence, the \textbf{GenInfCuts} introduces the feasibility cuts into the reduced master problem. If there doesn't exist any intra-layer cycle, then \textbf{GenUnavCuts} produces cuts/cutting-planes that alter the cost of the reduced master problem of future iterations to incorporate the effect of unavailabilities. Both variations of the proposed algorithm iteratively add the \textit{cutting-planes or cuts} (inequality constraints) to the reduced master problem, until the non-decreasing lower bound converges with the non-increasing upper bound. 


\subsubsection*{The reduced master problems}

The $n^{th}$ iteration of the reduced master problems for the two proposed approaches are described in \ref{form:rmp1} and \ref{form:rmp2} respectively. 
Both \ref{form:rmp1} and \ref{form:rmp2} comprise of spatial-only constraints from \ref{form:rpptu3if} i.e. these constraints are independent of temporal variables $\Gamma$; see $\mathcal{S}_X$ in Table \ref{tab:formnotas} and Equations \eqref{eqn:rmp:bounds} and \eqref{eqn:rmp:bounds2}. Additionally, \ref{form:rmp1} and \ref{form:rmp2} comprise of constraints expressed using $\gamma$, $d_k$, $\lambda$, and the degree constraints, see Equations \eqref{eqn:rmp:costtp}, \eqref{eqn:rmp:unav}, \eqref{eqn:rmp:hatlam}, and \eqref{eqn:rmp:deg}. 
\setcounter{eqc}{\theequation}
\stepcounter{eqc}
\begin{subequations}
\makeatletter
\def\@currentlabel{RP$^{(n)}$-\arabic{eqc}}
\makeatother
\label{form:rmp1}
\renewcommand{\theequation}{RP$^{(n)}$-\arabic{eqc}\alph{equation}}
 \begin{alignat}{3}
	\label{eqn:rmp:cost}
	&& \min_{\gamma, X, d_1, \dots, d_\mathcal{|\mathcal{K}|}} & \sum_{m,k,l} c_{m,k,l} X_{m,k,l} + \beta \gamma  && \\
	\label{eqn:rmp:costtp}
	&& \text{s.t.} ~~~ & \gamma \geq \sum_{m,l} t_{m,k,l} X_{m,k,l} + d_k, && \forall k \in \mathcal{K} \\ 
	\label{eqn:rmp:unav}
	& & & d_k \geq \sum_{m,l} \lambda_{m,k,l} X_{m,k,l}, &&  \forall \lambda \in \tilde{\mathcal{S}}_k^{(n-1)}, k \in \mathcal{K} \\ 
	\label{eqn:rmp:hatlam}
	& & & \sum_{m,k,l} \lambda_{m,k,l} X_{m,k,l} \leq \sum_{m,k,l} \lambda_{m,k,l} - 1,  ~~~ &&  \forall \lambda \in \hat{\mathcal{S}}_k^{(n-1)} \\ 
	\label{eqn:rmp:bounds}
	& & & X \in \mathcal{S}_X \textit{ and } d_k \geq 0 && 
 \end{alignat}
\end{subequations} 

The lower bound for RPP-TRU is contributed by the spatial cost and a temporal cost estimate $\gamma$ in Equations \eqref{eqn:rmp:cost} and \eqref{eqn:rmp2:cost}. The temporal cost estimate $\gamma$ accounts for the maximum running and waiting/delay time among the agents, see Equation \eqref{eqn:rmp:costtp}. In particular, the left hand side is the completion time of agent $k$, computed as the sum of the total running-time $\sum_{m,l} t_{m,k,l} X_{m,k,l}$ and the total delay experienced $d_k$. These $|\mathcal{K}|$ inequalities along with minimization over $\gamma$ ensure that $\gamma$ attains the maximum completion time among all the agents. The waiting/delay is achieved through the variables $d_k$ in Equation \eqref{eqn:rmp:unav}. These inequalities in Equation \eqref{eqn:rmp:unav} represent a family of constraints called unavailability cuts. The delays experienced by spatial solutions is implemented through the set $\tilde{\mathcal{S}}_k^{(n-1)}$ in Equation \eqref{eqn:rmp:unav} that grow in every iteration. The infeasible spatial solutions are eliminated by the feasibility cuts implemented through the set $\hat{\mathcal{S}}_k^{(n-1)}$, in Equation \eqref{eqn:rmp:hatlam}. Both the sets 
$\tilde{\mathcal{S}}_k^{(0)}$ and $\hat{\mathcal{S}}_k^{(0)}$ are 
initialized as empty for the first iteration $n=1$. The elements $\lambda$ from the sets $\tilde{\mathcal{S}}_k^{(0)}$ and $\hat{\mathcal{S}}_k^{(0)}$ are of the same size as $X$. Lastly, the spatial solutions must belong to the set $\mathcal{S}_X$, see Table \ref{tab:formnotas}. This introduces spatial constraints to the \ref{form:rmp1} from the formulation \ref{form:rpptu3if}. 
\setcounter{eqc}{\theequation}
\stepcounter{eqc}
\begin{subequations}
\makeatletter
\def\@currentlabel{RP$^{(n)}$-\arabic{eqc}}
\makeatother
\label{form:rmp2}
\renewcommand{\theequation}{RP$^{(n)}$-\arabic{eqc}\alph{equation}}
 \begin{alignat}{3}
	\label{eqn:rmp2:cost}
	&& \min_{\gamma, X, d_1, \dots, d_\mathcal{|\mathcal{K}|}} & \sum_{m,k,l} c_{m,k,l} X_{m,k,l} + \beta \gamma  && \\
	\nonumber 
	&& \text{s.t.} ~~~ & \textit{Equations \eqref{eqn:rmp:costtp} and \eqref{eqn:rmp:unav} } \\ 
	\label{eqn:rmp:deg} 
	& & & \sum_{a_{m,k,l'} \in \delta^+_{ikl}} X_{mkl'} \leq 1; && \forall v_{ikl} \in \mathcal{V} \\ 
	\label{eqn:rmp:bounds2}
	& & & X \in \mathcal{S}_X \textit{ and } d_k \geq 0 && 
 \end{alignat}
\end{subequations} 

In the reduced master problem formulation \ref{form:rmp2}, Equation \eqref{eqn:rmp:deg} represents degree constraints, that restricts more than one incoming arc with integer decision. This modification eliminates need of the feasibility cuts given by Equation \eqref{eqn:rmp:hatlam}. The necessity of the feasibility cuts in \ref{form:rmp1} is established in the literature \cite{conejo}, while its redundancy in \ref{form:rmp2} is established later in Section \ref{rpptu:sec:analysis} through Remark \ref{rem:dominance}. 

In Section \ref{rpptu:sec:analysis}, we show that the feasible region of \ref{form:rmp1} and \ref{form:rmp2} shrinks with increase in $n$; and eventually for some $n=N$, all the formulations \ref{form:rmp1}, \ref{form:rmp2}, and \ref{form:rpptu3if} have the same feasible region. The cost expressions of \ref{form:rmp1} and \ref{form:rmp2} are lower bound to the cost expression of the RPP-TRU, denoted as \ref{form:rpptu3if}, for all $n$ and spatial solution $X$, eventually becoming same at $n=N$. Briefly, this convergence of the reduced master problems to \ref{form:rpptu3if} is achieved by iteratively accounting for the temporal constraints of \ref{form:rpptu3if} from Equation \eqref{form:run} and the temporal cost contribution $\beta \gamma$ from Equation \eqref{form:rpptucost} using the feasibility cuts from \textbf{GenInfCuts} block and the unavailability cuts from \textbf{GenUnavCuts} block respectively. 

The \textbf{reduced master problems} \ref{form:rmp1} and \ref{form:rmp2} are indicated as blocks in Figures \ref{fig:bndrs1_rpp-tu} and \ref{fig:bndrs2_rpp-tu} respectively. These problems may be solved using any standard MILP solver, column generation approach, branch-and-bound approach, etc. In this work, we use Matlab's MILP solver to solve \ref{form:rmp1} and \ref{form:rmp2} because plenty of commercial and easy-to-use solvers are available in the market, which are specifically suited for RPP-TRU over small-sized networks. 

\subsubsection*{The temporal sub-problems} 
\label{rpptu:sec:tmpsp} 

The temporal sub-problems \ref{form:bsp}, for each agent $k \in \mathcal{K}$, is a consequence of decomposing \ref{form:rpptu3if} into a spatial part and multiple temporal parts. 
The temporal portion of \ref{form:rpptu3if} is comprised of running-time constraints in Equation \eqref{form:run} and unavailability constraints in Equation \eqref{form:unav}, that are decomposed into the temporal sub-problems. In particular, there are $|\mathcal{K}|$ sub-problems for each agent $k$ in each iteration $n$ that solve the temporal part of RPP-TRU given by Equation \eqref{form:run}, provided the spatial solution is a known value denoted as $\bar{X}^{(n)}$. These decoupled sub-problems are solved separately. 
\setcounter{eqc}{\theequation}
\stepcounter{eqc}
\begin{subequations}
\makeatletter
\def\@currentlabel{TSP$_k^{(n)}$}
\makeatother
\label{form:bsp}
\renewcommand{\theequation}{TSP$_k^{(n)}$-\arabic{eqc}\alph{equation}}
 \begin{alignat}{2}
	\label{form:bsp:cost}
 	\min_{\Gamma} ~&\Gamma_{dk0} \\ 
	\label{form:bsp:run} 
	\text{s.t.} ~& \Gamma_{jkl'} \geq t_{mkl} + \Gamma_{ikl}; && \forall \bar{X}^{(n)}_{mkl} = 1 \text{ in the $n^{th}$ iteration} \\ 
	\label{form:bsp:unav}
 	& \left.
    \begin{array}{ll}
     either,&\\
     ~~\Gamma_{jkl} \leq &z^- - t_{mkl} \\
     or, &\\
     ~~\Gamma_{jkl} \geq &z^+  \\
    \end{array} 
    \right\} && \forall (z^-, z^+) \in Z, \text{ if $\bar{X}^{(n)}_{mkl} = 1$} \\ 
	\label{form:bsp:intlin}
 	& \Gamma_{jkl} \geq 0 
 \end{alignat}
\end{subequations} 
Equation \eqref{form:bsp:run} is derived from Equation \eqref{form:run}, and it is only imposed if the arc $a_{mkl}$ connecting the head vertex $v_{jkl'}$ and tail vertex $v_{ikl}$ has decision $\bar{X}^{(n)}_{mkl} = 1$. Similarly, Equation \eqref{form:bsp:unav} is derived from Equation \eqref{form:unav}, and it is only imposed if the arc $a_{mkl}$, with tail vertex $v_{jkl}$, has decision $\bar{X}^{(n)}_{mkl} = 1$. Note that, in the first approach \ref{form:rmp1}, the spatial solution $\bar{X}^{(n)}$ of the $n^{th}$ iteration might make \ref{form:bsp} infeasible. 
In the second approach \ref{form:rmp2}, the spatial solution always yields a feasible temporal solution i.e. \ref{form:bsp} always has a finite optimizer. 

\ref{form:bsp} may be solved as a linear programming problem by applying dummy variables as shown in Equation \eqref{eqn:unav_dummy}. However, for computational efficiency, we apply a Depth First Search (DFS) like approach to determine the temporal values $\bar{\Gamma}^{(n)}$ of the $n^{th}$ iteration. The proposed DFS-like approach is suitable for determining delayed paths and generating new unavailability cutting-planes for \ref{form:rmp1}. In addition, this approach eliminates $\tau$ from the running-time and unavailability constraints in Equation \eqref{form:run} and \eqref{form:unav}.  
The temporal solution $\bar{\Gamma}^{(n)}$ is derived from Algorithm \ref{algo:temporal} and an upper bound to the optimal solution of RPP-TRU is achieved as $\sum_{m,k,l} c_{mkl} X_{mkl} + \beta \max_{k \in \mathcal{K}} \Gamma_{dk0}$. This algorithm is illustrated as \textbf{TemporalSol} block in Figures \ref{fig:bndrs1_rpp-tu} and \ref{fig:bndrs2_rpp-tu}. Note that, Algorithm \ref{algo:temporal} is executed over a fully connected graph $\mathcal{G}_{\bar{X}^{(n)}}$, which has a maximum degree of $1$ (at most one incoming/outgoing arc per vertex) for feasible spatial solution $\bar{X}^{(n)}$.

\ifsubbuild
\begin{algorithm}[H]
\else
\begin{algorithm}
\fi
	\caption{\textbf{TemporalSol}; Temporal solution for agent $k \in \mathcal{K}$ and Benders' iteration $n$} 
	\label{algo:temporal}
	
	\SetKwInOut{Input}{input}
	\SetKwInOut{InOut}{}
	\SetKwInOut{Output}{output}
	\Input{Solution sub-graph $\mathcal{G}_{\bar{X}^{(n)}}$ (all arcs $a_{mkl}$ with decision $\bar{X}^{(n)}_{mkl}=0$ are removed from graph $\mathcal{G}$),} 
	\InOut{Travel time $t_{mkl}$, and} 
	\InOut{Unavailability list for each arc $Z_m$} 
	\Output{Temporal solution $\bar{\Gamma}^{(n)}$ of $n^{th}$ iteration, and}
	\InOut{Arc and delay sequence, $aseq$ and $dseq$ respectively} 
	\BlankLine 
	$aseq, dseq \leftarrow \emptyset$ \tcp*{to store the sequence of arcs and corresponding delays} 
	$\Gamma_{ikl} \leftarrow 0, ~\forall v_{ikl} \in \mathcal{V}$ \tcp*{initialize the temporal solution} 
	$cv \leftarrow v_{sk0}$ \tcp*{the current vertex of $k^{th}$ agent is set as the source vertex} 
	$st \leftarrow 0$ \tcp*{start time initialization} 
	\While{true}{
		$a_{mkl} \leftarrow \delta^-(cv)$ $~~~$ (using graph $\mathcal{G}_{\bar{X}^{(n)}}$) \tcp*{arc outgoing at the current vertex is determined} 
		$rt \leftarrow t_{mkl}$ \tcp*{the `\emph{running-time}' of the arc $a_{mkl}$ is extracted } 
		\eIf{$(st, rt)$ is violated by the $\sigma^{th}$ entry in $Z_m$}{
		    \tcp{see Figure \ref{fig:unav} for checking violation }
		    $(z^-,z^+) \leftarrow Z_m(\sigma)$ 
		}{
                $z^+ \leftarrow st$
            }
		$dt \leftarrow z^+ - st$ \tcp*{`\emph{delay time}' (waiting) experienced at $cv$ } 
		$cv \leftarrow F^+(a_{mkl})$ \tcp*{start vertex is updated as the head vertex of arc $a_{mkl}$} 
		$st \leftarrow z^+ + rt$ \tcp*{starting time for the new current vertex $cv = F^+(a_{mkl})$} 
		$aseq \leftarrow aseq ~\cup~ a_{mkl}$ \tcp*{store the current arc} 
		$dseq \leftarrow dseq ~\cup~ dt$ \tcp*{store the corresponding delay} 
		$\Gamma_{ikl} \leftarrow st$ $~~~$ (here $v_{ikl} = cv = F^+(a_{mkl})$) \tcp*{the arrival time at the head vertices of each arc in the spatial solution is recorded} 
		\If{$cv$ is the terminal vertex $v_{dk0}$}{
		    break the loop 
		} 
	}
\end{algorithm}

\textbf{The unavailability cutting-planes:} The unavailability cuts are determined from Algorithm \ref{algo:unav} using the delay information obtained from Algorithm \ref{algo:temporal}. This algorithm is illustrated as \textbf{GenUnavCuts} in Figures \ref{fig:bndrs1_rpp-tu} and \ref{fig:bndrs2_rpp-tu}. 
Equation \eqref{eqn:unavcut2} shows the mathematical structure of the unavailability cutting-planes. 
\begin{equation}
	\label{eqn:unavcut2}
	\begin{aligned}
		d_k & \geq \Delta \Gamma \big(X_{mkl} - \sum_{a_{m'kl'} \in \delta^-(P_{ikl})} X_{mkl} \big) \ \forall k \in \mathcal{K}\\
		\textit{and } d_k & \geq 0 \ \forall k \in \mathcal{K}
	\end{aligned}
\end{equation}

The coefficients of spatial variables $X$ from Equation (\ref{eqn:unavcut2}) are stored in $\tilde{\mathcal{S}}_k^{(n)}$ for every feasible Benders' iteration. 
Here, $P_{ikl}$ is a path from the start vertex $v_{sk0}$ to vertex $v_{ikl}$, having unavailability violated. The set $\delta^-(P_{ikl})$ represents the set of arcs incoming at the path, but excluding the arcs involved in the path itself. The arc $a_{mkl}$ is the outgoing arc at the vertex $v_{ikl}$ (i.e. $v_{ikl} = F^-(a_{mkl})$) and having $X_{mkl}=1$. The delay $\Delta \Gamma$ is computed based on the cumulative delay caused at the head vertex of arc $X_{mkl}$ i.e. using the total delay in the path $P_{ikl} \cup a_{mkl}$. The constraint is visualized in Figure \ref{fig:unav_viz}. Observe that, only the spatial solutions that have the path $P_{ikl}$ from $v_{sk0}$ to $v_{ikl}$ will have $\Delta \Gamma \big(X_{mkl} - \sum_{a_{m'kl'} \in \delta^-(P_{ikl})} X_{mkl} \big) = \Delta \Gamma$. Any other path will be incoming through at least one of the arcs in $\delta^-(P_{ikl})$, resulting in $\Delta \Gamma \big(X_{mkl} - \sum_{a_{m'kl'} \in \delta^-(P_{ikl})} X_{mkl} \big) \leq 0$. Lastly, if all agents have the same running-time data, unavailability cuts for other agents $k' (\neq k) \in \mathcal{K}$ can be generated using the cut for agent $k$ given in Equation \eqref{eqn:unavcut2}. 

\begin{figure}
    \centering
    \includegraphics[width=0.8\textwidth]{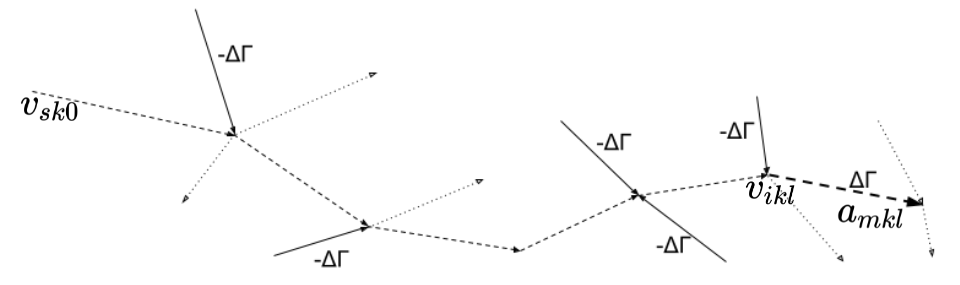}
    \caption{Visualizing the coefficients of the unavailability cutting-planes: The dashed arrows indicate the arcs from a portion of the solution $\bar{X}^{(n)}$ of the $n^{th}$ iteration, denoted by path $P_{ikl} \cup \{a_{mkl}\}$. The solid arrows indicate the arcs incoming at the vertices of the path $P_{ikl}$. All other arcs are denoted with dotted arrows. Since delay was experienced in the arc $a_{mkl}$, all incoming arcs to the path $P_{ikl}$, in between vertex $v_{sk0}$ and vertex $v_{ikl}$, denoted by solid arrows and given by the set $\delta^-(P_{ikl})$, have negative $\Delta \Gamma$ coefficients. Only arc $a_{mkl}$ has a positive $\Delta \Gamma$ coefficient. }
    \label{fig:unav_viz}
\end{figure}

Note that, the solution set $S_X$ will contain one cycle connecting the layers and depot (inter-layer cycle), and optionally one or more intra-layer cycles (contained entirely in one layer) that is either connected/disconnected to the inter-layer cycle. Observe that, Remark \ref{rem:unav_valid} insists that the connected intra-layer cycles must not be allowed in the solution. This is taken care of by either the degree constraints or the feasibility cuts. By design the unavailability cuts raise the lower bound to account for the delay in the path. However, if a cycle is connected to this path, the unavailability cut is no longer active i.e. the effect of delay in the path is not added to $d_k$ and therefore the cost of the path is incorrect in future iterations. 

\begin{rem}
    \label{rem:unav_valid}
    The unavailability cuts are effective only if the sub-graph constructed from the solutions of the reduced master problem doesn't have connected intra-layer cycles. 
\end{rem}

\ifsubbuild
\begin{algorithm}[H]
\else
\begin{algorithm}
\fi
	\caption{\textbf{GenUnavCuts}; Generate unavailability cuts for each agent $k \in \mathcal{K}$ and Benders' iteration $n$} 
	\label{algo:unav}
	\SetKwInOut{Input}{input}
	\SetKwInOut{InOut}{}
	\SetKwInOut{Output}{output}
	\Input{Replicated graph $\mathcal{G}$, }  
	\InOut{Arc and delay sequence, $aseq$ and $dseq$ respectively, and} 
	\InOut{Set of previous unavailability cuts $\tilde{\mathcal{S}}_{k,n-1}$}
	\Output{Updated set of unavailability cuts $\tilde{\mathcal{S}}_k^{(n)}$} 
	\BlankLine 
	\For{each arc $a_{mkl}$ in $aseq$}{
	    $dt \leftarrow dseq$ corresponding to $a_{mkl}$ \tcp*{extract delay corresponding to arc $a_{mkl}$} 
		\If{$dt$ is zero}{
		    continue to the start, to select next arc 
		} 
		$\Delta \Gamma \leftarrow$ cumulative sum of $dseq$ till arc $a_{mkl}$ \tcp*{sums up total delay experienced till the current arc $a_{mkl}$ in $aseq$} 
		$v_{ikl} \leftarrow F^-(a_{mkl})$ \tcp*{find tail vertex of the arc $a_{mkl}$} 
		Determine the set of arcs $\delta^-(P_{ikl})$ using $\mathcal{G}$ \tcp*{all incoming arcs to the tail vertices of sequence $aseq$ till vertex $v_{ikl} = F^-(a_{mkl})$, excluding the arcs of $aseq$ itself} 
		Initialize $\lambda$, i.e. $\lambda_{mkl} \leftarrow 0, ~\forall a_{mkl} \in \mathcal{A}$ \tcp*{initialize the cut with zeros} 
		$\lambda_{mkl} \leftarrow \Delta \Gamma$ \tcp*{set the cumulative delay as coefficient for current arc $a_{mkl}$} 
		$\lambda_{m'kl'} \leftarrow -\Delta \Gamma, ~\forall a_{m'kl'} \in \delta^-(P_{ikl})$ \tcp*{set the negative of cumulative delay as coefficient} 
		$\tilde{\mathcal{S}}_k^{(n)} \leftarrow \tilde{\mathcal{S}}_{k,n-1} ~\cup~ \lambda$ \tcp*{add the cut to the set of unavailability cuts for agent $k$} 
	    \BlankLine 
	    \BlankLine 
		\tcp{Each cut is valid for all $k' \in \mathcal{K}$ in case of homogeneous agents} 
		\For{$k' \in \mathcal{K}$ and $k' \neq k$}{ 
    		Initialize $\lambda$, i.e. $\lambda_{mkl} \leftarrow 0, ~\forall a_{mkl} \in \mathcal{A}$ \tcp*{initialize the cut with zeros} 
    		$\lambda_{mk'l} \leftarrow \lambda_{mkl}, ~\forall a_m \in A ~\&~ l \in \mathcal{L}$ \tcp*{copy the cut for agent $k'$ } 
    		$\tilde{\mathcal{S}}_{k',n} \leftarrow \tilde{\mathcal{S}}_{k',n-1} ~\cup~ \lambda$ \tcp*{add the cut to the iteratively growing set of unavailability cuts for agent $k'$} 
    	}
	}
\end{algorithm}

\textbf{The feasibility cutting-planes:} 
Equation (\ref{eqn:rmp:hatlam}) represents the feasibility cutting-planes. If the spatial solution of the $n^{th}$ iteration of \ref{form:rmp1} has intra-layer cycles (cycles contained in any one layer of replicated graph), then for each cycle one infeasibility cut coefficient vector $\lambda$ is constructed. 
In particular, if there are cycles in the solution sub-graph that are contained within any one layer of the sub-graph, then $\lambda$ is constructed for each of these intra-layer cycles. For a particular cycle in layer $l$, the value of $\lambda_{mkl}$ is set as $1$ provided the value $\bar{X}^{(n)}_{mkl}$ is $1$. 
These vectors $\lambda$ are stored in the sets $\hat{\mathcal{S}}_k^{(n)}$, based on the corresponding agent-sub-graph. 
The feasibility cutting-planes ensures that these cycles are not a part of the solution in future iterations of the reduced master problem. It is implemented as the \textbf{GenInfCuts} block in Figure \ref{fig:bndrs1_rpp-tu}.

\subsection{Theory and analysis of the modification process}
\label{rpptu:sec:analysis}

Figure \ref{fig:bndrs_mod} summarizes the relation between each step of the Benders' algorithm and the proposed modifications for solving RPP. Since the unavailability constraints alter only the cost and not the feasible region of RPP-TRU (see Remark \ref{rem:feasball}), the proposed algorithm is developed in two stages: (1) simplifying each step of Benders' algorithm for RPP, and then (2) for RPP-TRU, adding unavailability cutting-planes (inequality constraints) to accommodate the effect of delay on the cost, as shown in the bottom portion of Figure \ref{fig:bndrs_rpp-tu}. In this subsection, we describe the modification process through relevant theoretical results and proofs. 
The proposed decomposition algorithm 
is developed from the standard framework of the Benders' decomposition algorithm, see \cite{conejo}. In this subsection, we also present a comprehensive insight into the convergence mechanism. In particular, we show that the lower bound produced by the reduced master problem and the upper bounds produced on solving the temporal sub-problems eventually converge to the optimal value. 

\begin{figure}
    \centering
    \includegraphics[width=0.8\textwidth]{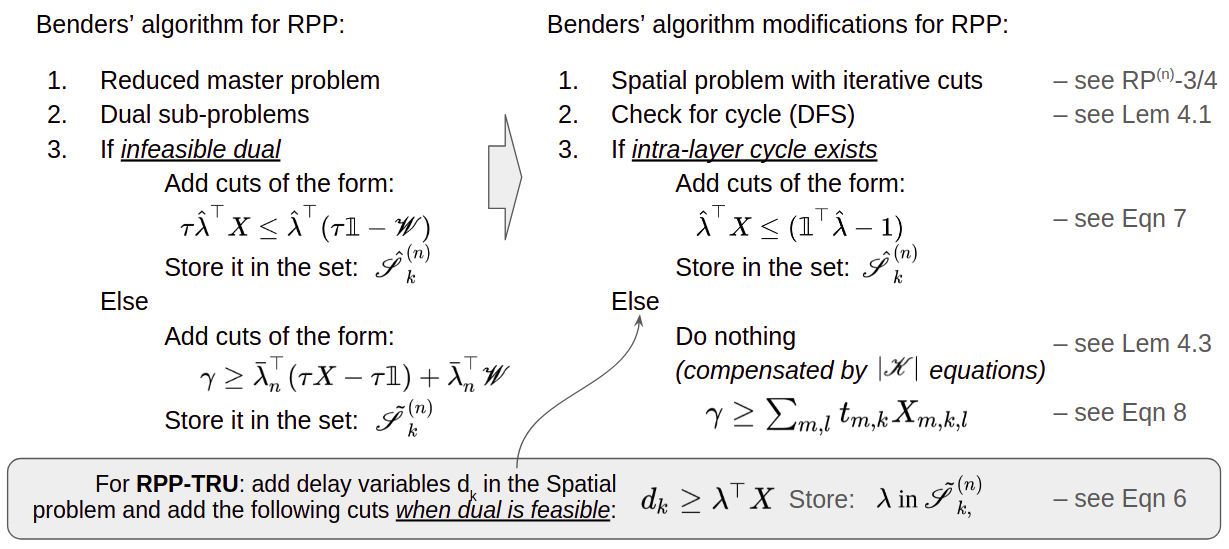}
    \captionsetup{justification=centering}
    \caption{Summary of modification of each step in Benders' algorithm for RPP (no temporal unavailabilities). In the proposed decomposition algorithm for RPP-TRU, delay variables $d_k$ are added, and the unavailability cuts are stored in $\tilde{\mathcal{S}}_k^{(n)}$. }
    \label{fig:bndrs_mod}
\end{figure}


\subsubsection*{Modification of the standard Benders' algorithm} 


The standard Benders' decomposition algorithm has three main steps as illustrated in Figures \ref{fig:bndrs_rpp-tu} and \ref{fig:bndrs_mod}. These steps are: (1) solving the reduced master problem, (2) solving the dual sub-problem, and lastly, (3) determining cuts for the reduced master problem of the next iteration based on feasible or infeasible dual solution. The reduced master problem for RPP-TRU, shown in \ref{form:rmp1} and \ref{form:rmp2}, is obtained by adding delay variables to the following reduced master problem for RPP: 
\begin{equation}
    \nonumber 
    \begin{aligned}
        \min_{\gamma, X} &~~\sum_{m,k,l} c_{m,k} X_{m,k,l} + \beta \gamma \\
        \text{   s.t. } &~~\text{Benders' cuts, expressed in terms of $\gamma$ and $X$ variables, are satisfied (see Figure \ref{fig:bndrs_mod}). } \\
        &~~X \in S_X  
    \end{aligned}
\end{equation}

The dual sub-problem for each agent $k \in \mathcal{K}$ is achieved from the temporal sub-problem (see \ref{form:bsp} in Section \ref{rpptu:sec:tmpsp}) for a given spatial solution $\bar{X}^{(n)}$ of the $n^{th}$ iteration. Lemma \ref{lem:fcut} shows that \textbf{infeasibility of the dual sub-problem} is determined by the existence of cycles in the layers of the solution sub-graph. Hence, a simple search algorithm is suitable for determining the feasibility of the dual sub-problem.   

\begin{lem} 
    \label{lem:fcut}
    In case of an empty unavailability list, each unbounded extreme ray in the dual sub-problem corresponds to a unique cycle 
    in the solution subgraph given by  $\bar{X}^{(n)} \in \mathcal{S}_X$, and vice versa. 
\end{lem}

\begin{prf} 
    In case of an empty unavailability list (for RPP), the constraints of the primal sub-problem are comprised only of the running-time constraints in Equation \eqref{form:run}. 
    The running-time constraints is equivalently expressed as $\mathcal{B}^\top \Gamma \geq \mathcal{W} - \tau (\mathbbm{1} - X)$; where $\mathcal{B}$ is the incidence matrix and $\mathcal{W}$ is the vector of running-time values of all arcs of the replicated graph $\mathcal{G}$ without the inter-component arcs (arcs connecting the agent-sub-graphs, see Figure \ref{fig:prelim_fig3} in Section \ref{prelims:sec:repg}). In particular, \textit{we assume without any loss of generality that} the replicated graph represented by the $\mathcal{B}$ matrix has $|K|$ disconnected components (agent-sub-graphs)\footnote{This assumption of disconnected agent-sub-graphs doesn't affect the solution. We make this assumption to keep the temporal flow constraints simple, i.e. no cycles allowed. \label{fn:assumpB}}. 
    
    The dual sub-problem, obtained from the feasibility temporal sub-problem i.e. running-time constraints with zero cost, reduces to the following linear program: 
    \begin{equation} \tag{$\textit{DP}(\bar{X}^{(n)})$}
	\label{prb:dsp}
	\nonumber
	\begin{aligned}
    	& \underset{\lambda}{\text{max}} & & (\tau \bar{X}^{(n)} - \tau \mathbbm{1} + \mathcal{W})^\top \lambda \\
    	& \text{s.t.} & & {\mathcal{B}} \lambda \leq 0 \\
    	& & & \lambda \geq 0
	\end{aligned}
    \end{equation} 
    Note that, if \ref{prb:dsp} is unbounded (maximum value is infinity) then, by the Linear Programming duality argument, there is no solution in the primal counterpart i.e. there doesn't exist a feasible temporal solution $\Gamma^{(n)}$ for given spatial solution $\bar{X}^{(n)}$. 
	
    The proof begins by establishing a relation between the dual vectors $\lambda$ and flow in the graph given by incidence matrix $\mathcal{B}$. 
    As per \ref{prb:dsp}, each entry in $\lambda$ correspond\authar[s] to one column of matrix $\mathcal{B}$ and hence one arc of graph $\mathcal{G}$. The value of row in vector $\lambda$ can be considered a flow value in the corresponding arc; thus the equation $\mathcal{B}\lambda = 0$ becomes a flow conservation constraint (similar to Equation (\ref{form:flow})). In flow constraints given by Equation \eqref{form:flow}, the product $\mathcal{B}\lambda$ is evaluated against $0$ vector, and hence the sum of incoming and outgoing flow values at any vertex must be zero. 
	
    {\authar[The `if' part of this] lemma is proved using the following claims: }
    \begin{itemize}
        \item \authar[Claim 1:] the sets $\{\lambda \geq 0 | \mathcal{B}\lambda \leq 0 \}$ and $\{\lambda \geq 0 | \mathcal{B}\lambda = 0 \}$ represent the same region, 
        \item \authar[Claim 2:] feasible region represented by the set $\{\lambda \geq 0 | \mathcal{B}\lambda = 0 \}$ is a cone, 
        \item \authar[Claim 3:] any $\lambda \in \{\lambda \geq 0 | \mathcal{B}\lambda \leq 0\}$, that does not correspond to a cycle, is not an extreme ray.
    \end{itemize}
	
	\authar[Claim 1] is evident, if \authar[the vector] $\mathcal{B}\lambda$ cannot have identically non-positive entry for any $\lambda \geq 0$. This statement can be argued by contradiction using two exhaustive cases: firstly, flow cannot originate at a vertex as this implies positive right hand side corresponding to that vertex (or row); and secondly, cyclic flow having imbalance implies there is a source and sink for the surplus value, making atleast one entry positive and one negative in the right hand side. 
	
	\authar[Claim 2] is true because any feasible vector $\lambda$ can be scaled by $\alpha \in [0,\infty)$ and \authar[it ]still remain\authar[s] feasible. As a consequence, the entire feasible region can be described using conic combination of extreme rays. 
	
	For \authar[claim 3], the decomposition property of flow is utilized. 
	Consider $\lambda \in \{\lambda \geq 0 | \mathcal{B}\lambda \leq 0 \}$, and construct a sub-graph having only those arcs that have non-zero entry in $\lambda$. 
	If the in-degree or out-degree for any vertex in this sub-graph is greater than $1$, then neglect\authar[ing] all except one direction, select\authar[s] a cycle of the graph (it cannot be acyclic as our first claim holds). The minimum flow value of this cycle (say $\lambda_1$) can be subtracted from all the flow values of this cycle, thus reducing flow of some arcs to zero and eliminating them in the process. Since both $\lambda$ and $\lambda_1$ are in null space of $\mathcal{B}$, the resulting flow (say $\lambda_2 := \lambda - \lambda_1$) will also belong to null space of $\mathcal{B}$, aside from being non-negative. Hence, any feasible $\lambda$, that doesn't correspond to a cycle, can be represented using sum of two other vectors from the cone, implying that its not an extreme ray. In other words sub-graph corresponding to an extreme ray must be a cycle. 
	
	The `only-if' part of the proof can be established by showing that a cycle \authar[is not only] non-decomposable, but also relate\authar[d] to a single valued positive flow that satisfies flow constraints ($\mathcal{B}\lambda = 0$). In other words, $\lambda$, corresponding to any cycle, cannot be represented as a sum of any two elements of the cone, and hence it is an extreme ray. 
	
	Note that if an extreme ray, denoted by $\lambda$, relates to a cycle in the solution subgraph given by $\bar{X}^{(n)}$, then the dual cost is unbounded in that direction i.e. the dual cost $(\tau \bar{X}^{(n)} - \tau \mathbbm{1} + \mathcal{W})^\top (\alpha \lambda) \rightarrow +\infty$, as $\alpha \rightarrow +\infty$. Such unbounded extreme rays are chosen as binary $0$-$1$ vectors denoted by $\hat{\lambda}$ in this work.
\end{prf}

\textbf{In case there exists an intra-layer cycle} in the solution sub-graph $\mathcal{G}_{\bar{X}^{(n)}}$ i.e. dual sub-problem is infeasible then a feasibility-cut of the form $\hat{\lambda}^\top \left(\tau X - \tau \mathbbm{1} + \mathcal{W}\right) \leq 0$ is introduced to the reduced master problem of the next iteration $n+1$, see \cite{conejo}. 
This inequality is violated because $\alpha \hat{\lambda}$ (for some $\alpha > 0$) is an unbounded extreme ray of the $n^{th}$ iteration, and hence the optimal cost of the dual sub-problem \ref{prb:dsp}, given as $\alpha \hat{\lambda}^\top \left(\tau \bar{X}^{(n)} - \tau \mathbbm{1} + \mathcal{W}\right)$, is unbounded or infinity at $\alpha \rightarrow \infty$. Hence, the inequality is violated by $\bar{X}^{(n)}$ i.e. $\hat{\lambda}^\top \left(\tau \bar{X}^{(n)} - \tau \mathbbm{1} + \mathcal{W}\right) > 0$. 

From Lemma \ref{lem:fcut}, an equivalent representation of feasibility-cut is derived as given in Equation \eqref{eqn:feasbcut}, Equation \eqref{eqn:rmp:hatlam} and Step-3 (if-condition) of Figure \ref{fig:bndrs_mod}. Here, $\hat{\lambda}$'s are unbounded extreme rays corresponding to intra-layer cycles in solution subgraph $\mathcal{G}_{\bar{X}^{(n)}}$. They are stored in an iteratively increasing feasibility-cut set, $\hat{\mathcal{S}}_k^{(n)} = \hat{\mathcal{S}}_k^{(n-1)}  \cup \{\textit{all }\hat{\lambda}\textit{ from cycles in } \mathcal{G}_{\bar{X}^{(n)}}\}$ 
initialized by null set, $\hat{\mathcal{S}}^{(0)} = \{\null\}$. 
Note that the feasibility cuts are cycle-elimination constraints, that prevent the formation of certain cycles in spatial solutions of future iterations. These constraints are more commonly known as sub-tour elimination constraints in the literature of Travelling Salesman Problem, see \cite{mtz}. 
\begin{equation}
	\label{eqn:feasbcut}
	\begin{aligned}
		\hat{\lambda}^\top X \leq (\mathbbm{1}^\top \hat{\lambda} - 1)
	\end{aligned}
\end{equation}


Based on Lemma \ref{lem:fcut}, three observations are made: (a) Remark \ref{rem:cyc} highlights the repetitive structure of feasibility cuts, as it is useful for extending the computed cuts to other parts of the replicated graph; (b) Corollary \ref{corr:hatsnempty} shows that the feasibility cuts will not be required in a non-temporal setting; and finally, (c) Remark \ref{rem:degnocyc} and \ref{rem:dominance} shows that degree constraints make the feasibility cuts redundant. The second approach for the proposed Benders' algorithm \ref{form:rmp2} utilizes the degree constraints in Equation \eqref{eqn:rmp:deg} instead of feasibility cuts in Equation \eqref{eqn:rmp:hatlam}. Also, Corollary \ref{corr:hatsnempty} and Remark \ref{rem:degnocyc} together give us a MILP formulation to solve the RPP using \ref{form:rmp1} or \ref{form:rmp2} with no additional iterative cutting-planes, as stated in Remark \ref{rem:rpp}. 

\begin{rem}
\label{rem:cyc}
All cycles, that don't include source and destination vertices, are restricted to only one layer of \authar[a] component \authar[in replicated graph $\mathcal{G}$] (see Figure \ref{fig:prelim_fig3} for illustration of a replicated graph). Also, each cycle has a counterpart in every layer of all components in $\mathcal{G}$ due to its repetitive structure.  
\end{rem}


\begin{corr}
\label{corr:hatsnempty} 
In case of empty unavailability list, set of feasibility-cut $\hat{\mathcal{S}}_k^{(n)}$ is always empty. 
\end{corr}

\begin{prf}
From the argument presented in the proof of Theorem \ref{thm:wd}, any cycle present in any layer of the replicated graph can be collapsed to form a vertex. \authar[Contrary to the method used in proof of Theorem \ref{thm:wd}, the vertex is assigned the minimum time-stamp] among all time-stamps in the cycle\authar[. This] yields a feasible solution of lower cost i.e. this vertex-disjoint solution dominates all solutions with such intra-layer cycles. \authar[Again, s]ince there are no possible inter-layer cycles (by Remark \ref{rem:cyc}) or inter-component cycles (assuming disconnected replicated graph\footref{fn:assumpB} on page \pageref{fn:assumpB} of Lemma \ref{lem:fcut}), the spatial solution of the reduced master problem for RPP will never have cycles. 
\end{prf}

\begin{rem}
\label{rem:degnocyc} 
If the degree constraints (see Equation \eqref{eqn:rmp:deg}) and binary constraints on variables $X \geq 0$ are satisfied, then there exists no solution $X$ whose subgraph $\mathcal{G}_X$ has connected intra-layer cycles. In particular, only the solution with disconnected intra-layer cycles is feasible.
\end{rem}

\begin{rem}
\label{rem:dominance} 
Since the solutions with connected intra-layer cycles are eliminated by the degree constraints, all solutions with intra-layer cycles in the reduced master problem \ref{form:rmp2} are dominated by the solution without these intra-layer cycles (i.e. additional cycles increase the cost value). This dominance property is valid in the presence of the unavailability cuts due to Remark \ref{rem:unav_valid}. 
\end{rem}

\begin{rem}
\label{rem:rpp}
The first iteration (at $n=1$) of the reduced master problems \ref{form:rmp1} and \ref{form:rmp2} solves the RPP. 
\end{rem}

Observe that, optimal temporal solution for $n^{th}$ iteration $\bar{\Gamma}^{(n)}$ is determined from the temporal sub-problem \ref{form:bsp} in mainly two steps (see Algorithm \ref{algo:temporal}): first considering \authar[no unavailability case and solving] for $\bar{\Gamma}^{(n)}$ using the \textit{running-time data}; 
then increasing/\textit{delaying} entries of $\bar{\Gamma}^{(n)}$ minimally to achieve feasibility with respect to the unavailability constraints in polynomial time (Theorem \ref{thm:poly}). 

\begin{thm}
\label{thm:poly}
The sub-problem \ref{form:bsp} is solvable in polynomial time. 
\end{thm}

\begin{prf}
Presenting a polynomial time algorithm for updating temporal solution $\bar{\Gamma}^{(n)}$ 
suffices this claim. Initially, $\bar{\Gamma}^{(n)}$ is computed using running-time constraints only i.e. without the unavailability list based correction. This \authar[claim] is shown using the disjunctive hyper-cubes, described by constraints in Equation (\ref{form:bsp:unav})
, where the either-part consists of upper bounds while the or-part consists of lower bounds on the temporal variables $\Gamma_{jkl}$. Since the sub-problem \ref{form:bsp} is a minimization problem, only the lower bound given in the or-part can be tight. 

Taking the base case to be the portion of trajectory leading from source vertex to the depot vertex; this is always feasible with respect to unavailability constraints for all agents $k \in \mathcal{K}$. Let a portion of trajectory be directed from $v_{sk0}$ to $v_{jkl'}$ via $v_{ikl}$ for any arbitrary choice of $k$, as shown in Figure \ref{fig:dhc}. Assume (as an inductive hypothesis) that the portion of trajectory from $v_{sk0}$ to $v_{ikl}$ (excluding $v_{ikl}$) is feasible for all $k \in \mathcal{K}$ \authar[with respect to] the unavailability constraints, and the cost of this portion of the trajectory is minimum. If $v_{ikl}$ is minimally increased by $\Delta \Gamma$ to make $\Gamma_{ikl}$ feasible, as shown in the bottom-half of Figure \ref{fig:dhc}, then the portion of trajectory from $v_{sk0}$ to $v_{ikl}$ (including vertex $v_{ikl}$) is feasible \authar[with respect to] unavailability constraints, and the cost is also minimum by inductive hypothesis. In short, correction of temporal variables (given the spatial solution) doesn't affect the temporal data in variables occurring earlier in the trajectory. 

\begin{figure}[!h]
	\centering
	\includegraphics[width=0.6\textwidth]{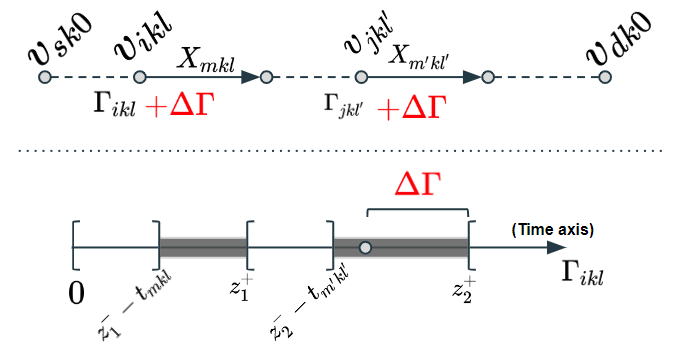}
	\captionsetup{justification=centering}
	\caption{The top figure shows a sample vertex disjoint solution trajectory for agent $k$. The bottom figure shows that the time value (marked with cicle) at vertex $\Gamma_{ikl}$ is infeasible \authar[because it lies in the unavailable region given by] $(z^-_2-t_{mkl}, z^+_2)$; \authar[where $t_{mkl}$ is the running-time of the outgoing occupied arc (see Figure \ref{fig:unav} to observe the effect of running-time on lower limit of unavailable region). H]ence minimum increment of $\Delta \Gamma$ is required to correct it. Since $\Gamma_{jkl'}$ appears later in this solution sequence, correction is also reflected in this variable.}
	\label{fig:dhc}
\end{figure}

Thus, this simple greedy algorithm (Algorithm \ref{algo:unav}) corrects the temporal variable $\Gamma$ in ascending order in $\mathbf{O}(|\mathcal{V}|)$ steps, as there are at most $|\mathcal{V}|$ temporal variables. 
There can only be a finite number of unavailabilities (say $\eta$) within the inspection time bound 
\authar[because] the availability gaps have finite non-zero lower bound. In particular, since temporal-graph agent tuple $(G,Z,W)$ is well-defined, every availability gap of an arc is at least as large as the running-time value for that arc. If the smallest running-time value is given as $t_{min}$ $(:= \min_{a_{mkl} \in \mathcal{A}} \{t_{mkl}\})$, and the number of availability gaps for each vertex $v_{ikl} \in \mathcal{V}$ is at \authar[most $\eta$, then] the product \authar[(given as $\eta t_{min}$)] must be upper bounded by the inspection time $T_I$; where Corollary \ref{cor:timebound} gives us $\Gamma_{ikl} \leq T_I, ~\forall v_{ikl} = \mathcal{V}$ and $T_I = (|V||\mathcal{L}|+1)(T_p + \max_{a_{mkl} \in \mathcal{A}} t_{mkl})$. Thus for each vertex in a given spatial solution, the number of comparisons with entries of unavailability list will be at most $\frac{T_I}{t_{min}}$. 
\end{prf}

\textbf{In case there doesn't exists any intra-layer cycle} in the solution subgraph $\mathcal{G}_{\bar{X}^{(n)}}$ of the $n^{th}$ iteration, then the dual sub-problem \ref{prb:dsp} has a feasible solution. Using Lemma \ref{lem:ccut}, we show in Equation \eqref{eqn:costcut} that the entire family of cost-cuts for RPP is replaced by a set of $|\mathcal{K}|$ cuts given as $\gamma \geq \sum_{m,l} t_{m,k,l} X_{m,k,l}$ for each $k \in \mathcal{K}$. In reduced master problem of RPP-TRU, the delay variables $d_k$ are also added to these $|\mathcal{K}|$ cuts, as shown in Equation \eqref{eqn:rmp:costtp} and Step-3 (else-condition) of Figure \ref{fig:bndrs_mod}. The unavailability cuts using $d_k$ and $X$ are constructed based on a path in the current spatial solution $\bar{X}^{(n)}$ as proposed in Equation \eqref{eqn:unavcut2}, using Algorithm \ref{algo:unav}. The unavailability cuts ensure that in future iterations of the reduced master problems \ref{form:rmp1} and \ref{form:rmp2}, any spatial solution with these paths will also have a delay-based component in the lower bound. 

\begin{lem}
    \label{lem:ccut}
    In case of no unavailability, for single-agent scenario ($|\mathcal{K}| = 1$), if the dual sub-problem \ref{prb:dsp}, given $\bar{X}^{(n)} \in \mathcal{S}_X$, has a finite maximum then the dual maximizer $\bar{\lambda}^{(n)}$ equals $\bar{X}^{(n)}$. 
\end{lem}

\begin{prf}
    The finite maximizer search problem is derived from \ref{prb:dsp}, thus the `no unavailability' scenario changes the feasible set in \ref{prb:dsp} to $\{\lambda \geq 0 | \mathcal{B}\lambda \leq \mathcal{C}\}$; where the coefficient of $\mathcal{C}$ corresponding to vertex $v_{d,1,0}$ is $1$ and rest are $0$. Since it is assumed that this optimization problem has a finite maximum, by Lemma \ref{lem:fcut}, the solution subgraph doesn't contain cycles. Also, as discussed earlier in the proof of Lemma \ref{lem:fcut}, positive cost in the dual sub-problem is only possible if non-zero entries of $\lambda$ are related to arcs of solution subgraph. Clearly, the solution subgraph is a vertex-disjoint path and the entries of $\lambda$ corresponding to the vertex-disjoint path is non-zero. In addition, the constraint set $\mathcal{B}\lambda \leq \mathcal{C}$ restricts the value of $\lambda$, at destination vertex $d$, to be at most $(\mathcal{C})_{d,1,0}$. By backtracking the destination vertex along the vertex-disjoint paths, it \authar[is] found out that all non-zero entries of $\lambda$ \authar[are] at most $(\mathcal{C})_{d,1,0}$. 
	
    Observe that $\bar{X}^{(n)}$ is an element of the set $\mathcal{S}_X$, which implies $(\mathcal{B}\lambda)_{d,1,0} \leq (\mathcal{C})_{d,1,0} := 1$ for all sink vertices in graph $\mathcal{G}$, and $(\mathcal{B}\lambda)_{vkl} \leq 0$ for all other vertices. Hence the choice of maximizer $\bar{\lambda}^{(n)} = \bar{X}^{(n)}$ produces maximum cost by selecting all arcs in the solution subgraph. 
\end{prf}

The cost-cuts can be equivalently represented as $\gamma\geq \bar{\lambda}^{(n)\top} (\tau X - \tau \mathbbm{1} + \mathcal{W})$, where $\bar{\lambda}^{(n)}$'s are picked from the iteratively increasing cost-cut set, $\bar{\mathcal{S}}_k^{(n)} = \bar{\mathcal{S}}_k^{(n-1)}  \cup \{\bar{\lambda}^{(n)}\}$; ($\bar{\mathcal{S}}^{(0)} = \{\phi\}$). From the standard framework for cost-cut in Benders' algorithm, the following is observed (provided $\tau$ is very large, see Corollary \ref{cor:timebound}):
\begin{equation}
    \nonumber 
    \begin{aligned}
        & & & \{\gamma\geq \bar{\lambda}^{(n)\top} (\tau X - \tau \mathbbm{1} + \mathcal{W}), \gamma\geq 0\} \\ 
        & \iff & & \{\gamma\geq \bar{\lambda}^{(n)\top} (\tau X - \tau \mathbbm{1}) + \bar{\lambda}^{(n)\top} \mathcal{W}, \gamma\geq 0\}\\
        & \iff & & \{\gamma \geq \bar{\lambda}^{(n)\top} \mathcal{W}\textit{ if $X = \bar{X}^{(n)}$, } 0 \textit{ otherwise}\}
    \end{aligned}
\end{equation}

Consequently, given $\bar{\lambda}^{(n)} = \bar{X}^{(n)}$ from Lemma \ref{lem:ccut}, we get the following relation for each agent $k \in \mathcal{K}$: 
\begin{equation}
    \label{eqn:costcut}
    \begin{aligned}
        \gamma \geq \left\{ \begin{matrix} \bar{X}^{(n)\top} \mathcal{W} &\textit{ if $X = \bar{X}^{(n)}$, }\\ 0 &\textit{ otherwise} \end{matrix} \right. 
    \end{aligned}
\end{equation}

Note that the feasibility-cut in Equations (\ref{eqn:feasbcut}) and \eqref{eqn:rmp:hatlam} and the cost-cut in Equations (\ref{eqn:costcut}) and \eqref{eqn:rmp:costtp} are independent of $\tau$. 

In a brief summary, the convergence of the proposed method is guaranteed by showing non-trivial equivalence with Benders' decomposition framework. The cost-cuts in Equation \eqref{eqn:rmp:costtp} are derived from Equation \eqref{eqn:costcut}. The cost-cuts and proposed unavailability cuts in Equation \eqref{eqn:rmp:unav} contribute to the improvement of the lower bound. The feasibility cuts in Equation \eqref{eqn:rmp:hatlam} are derived from Equation \eqref{eqn:feasbcut}. Equation \eqref{eqn:rmp:deg} eliminates the necessity of the feasibility cuts from the second approach of the proposed algorithm.

\subsubsection*{Sketch of the convergence}

In this subsection, we present a comprehensive insight into the convergence of the proposed algorithm. In brief, the following points sketch the convergence mechanism: 
\begin{itemize}
    \item[1.] \textbf{$\textit{feasb}(\text{\ref{form:rpptu3if}})$ is a subset of $\textit{feasb}(\text{\ref{form:rmp1}})$ and $\textit{feasb}(\text{\ref{form:rmp2}})$\footnote{The degree constraints areignored. } at $n=1$:} Observe that \ref{form:rmp1} and \ref{form:rmp2} at $n=1$ have the sets $\hat{\mathcal{S}}_k^{(0)}$ and $\tilde{\mathcal{S}}_k^{(0)}$ initialized as empty. Since $\gamma$ and $d_k$ are unbounded from above, only $X \in S_X$ governs the feasible region at iteration $n=1$. Hence, the feasible region of RPP-TRU is a subset of the reduced master problem \ref{form:rmp1} i.e. $\textit{feasb}(\text{\ref{form:rpptu3if}}) \subseteq \textit{feasb}(\text{\ref{form:rmp1}})$. In particular, the feasible region of RPP-TRU has additional constraints other than $\mathcal{S}_X$, making it a subset of the reduced master problem \ref{form:rmp1}. The same reasoning implies $\textit{feasb}(\text{\ref{form:rpptu3if}}) \subseteq \textit{feasb}(\text{\ref{form:rmp2}})$ at $n=1$, because \ref{form:rmp1} is same as \ref{form:rmp2} when the degree constraints are ignored. 
    \item[2.] \textbf{The cost of reduced master problems are at most the cost of RPP-TRU:} Lemma \ref{lem:ccut} shows that, without the unavailability constraints, the cost value of the reduced master problems are same as the cost value of RPP-TRU for all $X$. Therefore, on comparing Equations \eqref{eqn:rmp:cost} and \eqref{eqn:rmp:costtp} with Equations \eqref{form:rpptucost} and \eqref{form:maxtime}, we get $\sum_{m,l} t_{m,k,l} X_{m,k,l} = \Gamma_{dk0}$ and $d_k = 0$, for all agents $k \in \mathcal{K}$ (the spatial portion of the cost is same $\sum_{m,k,l} c_{mkl} X_{mkl}$). The unavailability cutting-planes in Equation \eqref{eqn:unavcut2} ensure that the value of $d_k$ is exactly equal to the waiting due to unavailabilities $\Delta \Gamma$ for at least the spatial solutions of the previous iteration $\bar{X}^{(1)}, \dots, \bar{X}^{(n-1)}, \bar{X}^{(n)}$. Hence, for the unavailability constrained problem, we have $\sum_{m,l} t_{m,k,l} X_{m,k,l} + d_k \leq \Gamma_{dk0}$, for all agents $k \in \mathcal{K}$ and spatial solution $X$, resulting in the cost of reduced master problem be at most the cost of RPP-TRU. 
    \item[3.] \textbf{$\textit{feasb}(\text{\ref{form:rpptu3if}})$ is a subset of $\textit{feasb}(\text{\ref{form:rmp1}})$ and $\textit{feasb}(\text{\ref{form:rmp2}})$ for all $n$:} New constraints are always added to the reduced master problems of the previous iteration, hence the feasible region keeps shrinking. Therefore, the minimum value of the reduced master problems are non-decreasing with respect to the iterations $n$. Since, the cost of reduced master problem be at most the cost of RPP-TRU and the RPP-TRU is a subset of the reduced master problem, the minimum cost of the reduced master problems are at most ($\leq$) the minimum cost of the actual problem. Figure \ref{fig:bendersset} illustrates this relation in detail. 
    \begin{figure}
        \centering
        \includegraphics[width=0.6\textwidth]{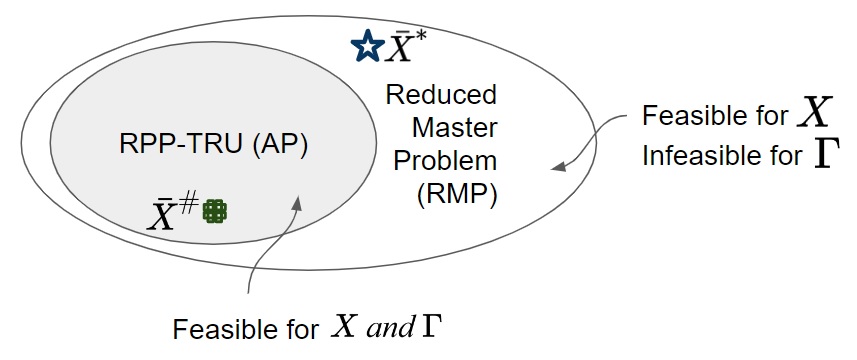}
        \caption{For brevity in illustration, we define the feasible set of the RPP-TRU as $\mathrm{P}$ and the reduced master problem as $\mathrm{RMP}$. Given: the feasible set of an actual problem $\mathrm{P}$ is a subset of a reduced master problem $\mathrm{RMP}$, the cost of actual problem $c_{\mathrm{P}}(X)$ is at least the cost of master problem $c_{\mathrm{RMP}}(X)$, and $\bar{X}^{\#}$ is the minimizer for $\mathrm{P}$. Case-1: If $\bar{X}^{\#}$ is also the minimizer for $\mathrm{RMP}$, then the costs satisfy the lower bound relation $c_{\mathrm{RMP}}(\bar{X}^{\#}) \leq c_{\mathrm{P}}(\bar{X}^{\#})$. Case-2: If $\bar{X}^*$ is the minimizer for $\mathrm{RMP}$, then the costs satisfy the lower bound relation $c_{\mathrm{RMP}}(\bar{X}^*) \leq c_{\mathrm{RMP}}(\bar{X}^{\#})$ provided $\bar{X}^{\#} \in \mathrm{RMP}$ for any arbitrary choice of $\bar{X}^{\#}$ i.e. $\mathrm{P} \subset \mathrm{RMP}$. Since $c_{\mathrm{RMP}}(\bar{X}^{\#}) \leq c_{\mathrm{P}}(\bar{X}^{\#})$, the optimal cost $c_{\mathrm{RMP}}(\bar{X}^*)$ is a lower bound of $c_{\mathrm{P}}(\bar{X}^{\#})$. }
        \label{fig:bendersset}
    \end{figure}
    \item[4.] \textbf{$\textit{feasb}(\text{\ref{form:rpptu3if}})$ is equal to $\textit{feasb}(\text{\ref{form:rmp1}})$ and $\textit{feasb}(\text{\ref{form:rmp2}})$ for some $n=N$:} The proposed algorithm is non-trivially equivalent to the Benders' decomposition algorithm, as shown in Section \ref{rpptu:sec:propbendalgo}. Therefore, the following conclusions are made: 
    \begin{itemize}
        \item[1.] The feasibility cuts are proxies to the running-time constraints. 
        \item[2.] The degree constraints eliminate the need of feasibility cuts. 
        \item[3.] The cost-cuts and unavailability cuts ensure that the temporal cost value of every element in the feasible region of the reduced master problem is same as that of the temporal cost value of RPP-TRU. In particular, the unavailability cuts are proxies to the temporal constraints and its effect on the cost value, while the cost-cuts update the cost value due to running-time. 
    \end{itemize} 
    The reduced master problems \ref{form:rmp1} and \ref{form:rmp2} have finite number of integer solutions because $X \in \{0,1\}^{|\mathcal{A}|}$. Therefore, there are a finite number of feasibility cuts and unavailability cuts, and at least one cut is added per iteration. Lemma \ref{lem:fcut} guarantees that by removing all the cycles from the solution of the reduced master problems, all infeasible spatial solutions are removed, and therefore the feasible region of the reduced master problems are same as the RPP-TRU in a finite number of steps. Similarly, the feasible solution costs are updated via the cost cuts and unavailability cuts in a finite number of steps. Consequently, there exist a finite $n=N$ for which, every spatial solution $X$ produces the same cost value for the reduced master problems and the RPP-TRU. 
\end{itemize} 

In conclusion the lower bounds are non-decreasing, and the reduced master problem shrinks to the RPP-TRU ensuring that the best upper bound is non-increasing. Eventually in a finite number of steps, the lower bound meets the upper bound, i.e. the optimal value of the reduced master problem is same as the optimal value of the RPP-TRU for the same spaital solution.

\section{Results and discussion}
\label{rpptu:sec:resdisc}

\ifsubbuild
\setcounter{section}{1}
\section{Results and discussion}
\fi




In this section, numerous simulation experiments are performed on randomly generated graph dataset. The details of the graph dataset are provided in the GitHub link: \cite{gitRPPTRU}. In addition, few results are also demonstrated on an example graph (see Figure \ref{fig:prelim_fig2}) first reported in \cite{lannez}. A part of Mumbai suburban railway network, Kurla-Thane-Vashi-Kurla (KTVK) network, is also modeled in a simulation case study for demonstrating real world application of the proposed approach. This section presents comparison results with suitable discussion on a couple of formulations, the algorithm approaches, and a real world application.

\subsection{Effect of time discretization} 

The comparison study shows us the benefits of the proposed mixed integer three-index formulation in contrast to time-discretization approach like TSN. 
In addition, the results also show that the proposed formulation size is independent of the planning time limit/bound. The results are demonstrated on a three-vertex-six-arc network picked from \cite{lannez} to provide a comparison analysis with sufficient details of TSN approach. Comparison over KTVK network shows the impact of unavailabilities on the proposed methodology and the TSN formulation. Lastly, a set of randomly generated graphs are used to show the approximate size comparison of the two formulation methodologies.  

\subsubsection*{Three-vertex-six-arc graph example}

TSN approach is implemented by discretizing time, and expanding the temporal graph to the time dimension, thus making it a static graph in the time-space dimension. The approach is detailed in \cite{pengThesis} and \cite{lannez}. This TSN approach is also called the static expansion of temporal graph \cite{othon}. 
A TSN based decision graph is constructed from the temporal graph-agent tuple $(G,Z,W)$, where $G$ represents the example graph with three vertices and six arcs, and $W$ represents the running-time data; as shown earlier in Figure \ref{fig:prelim_fig2}. The unavailabiltiy data $Z$ is shown in \textit{Ex-Ia} of Table \ref{tab:table1}. The TSN of this example network is illustrated in Figure \ref{fig:tsn3v6a}. The replicated graph for the three-vertex-six-arc graph is similar to the illustration shown in Figure \ref{fig:prelim_fig3}. Observe that for one agent and two agent scenario, as shown in \textit{Ex-I}, the replicated graph will only have one or two agent-sub-graphs respectively. 

\begin{figure}[!h]
	\centering
	\includegraphics[width=0.88\textwidth]{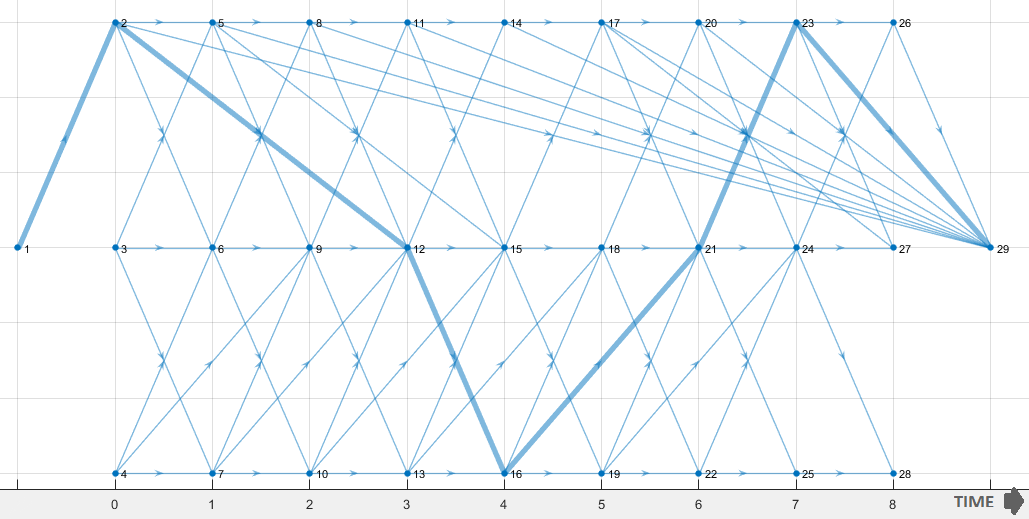}
	\captionsetup{justification=centering}
	\caption{TSN representation of a graph with three vertices and six arcs \authar[(see Figure \ref{fig:prelim_fig2})]. The x-axis represents time stamps discretized with interval of $1$ minute. The three vertices from Figure \ref{fig:prelim_fig2} are shown vertically (y-axis) for each time-stamp (x-axis). In particular copies of vertices $v_1, v_2, v_3$ are numbered as $3 t + 2, 3 t + 3, 3 t + 4$ respectively, where $t = \{0, \dots, 8\}$. These vertices are connected by arcs from left to right based on the running-time of that arc. An arc with running-time $3$ minutes connects vertices at $t=t_0$ and $t=t_0+3$; e.g. arc $a_2$ from Figure \ref{fig:prelim_fig2} with running-time $3$ minutes connects vertex marked as $2$ at $t=0$ to vertex marked as $12$ at $t=3$. The \textit{waiting arcs} connect same vertex of consecutive time-stamps, thus displayed by horizontally drawn arcs. Based on unavailability list (see Table \ref{tab:table1}, Ex-I, Case I), some arcs are not present at certain times. Arcs to the sink vertex (rightmost vertex, numbered $29$) are introduced for every time-stamp such that an agent at vertex $v_1$ can end the trajectory anytime, and thereby not be penalized for waiting. Finally, the optimal trajectory for Case I is highlighted. (The illustration is generated using Matlab)}
	\label{fig:tsn3v6a}
\end{figure}

Unavailability list for different cases is shown in the second column of Table \ref{tab:table1}, along with suitable comparison results. The unavailability list for \textit{Ex-Ia} is directly picked from the corresponding article by \cite{lannez}, while unavailabilty list for \textit{Ex-Ib} and \textit{Ex-Ic} is constructed to show different trajectory possibilities. In an optimal solution of \textit{Ex-Ia}, the agent visits service arc $a_2$ first (connects $v_1$ to $v_2$, numbered $2$ and $12$ in Figure \ref{fig:tsn3v6a} respectively), followed by arc $a_5$ (connects $v_3$ to $v_2$, numbered $16$ and $21$ in Figure \ref{fig:tsn3v6a} respectively); while the sequence is reversed for \textit{Ex-Ib} i.e. arc $a_5$ is visited first, followed by arc $a_2$. In \textit{Ex-Ic}, the optimal solution is improved due to the presence of an extra agent, i.e. the solver assigns one service arc to one agent. 

\begin{landscape}
    \begin{table}[!htp]
	\centering
	\captionsetup{justification=centering}
	\caption{Results comparing algorithm runtime, number of variables and constraints of TSN and the proposed method. \textit{Ex-I} is performed over the three-vertex-six-arc example graph, while \textit{Ex-II} is performed over the KTVK network. } 
	\label{tab:table1}
	\begin{center}		
		\begin{threeparttable}
		\begin{tabular}{c|ccc|c|c|cc|cc|cc}
			\hline
			& \multicolumn{3}{c|}{Unavailability list} & \multirow{2}{*}{Time limit } & \multirow{2}{*}{No. of agents} & \multicolumn{2}{c|}{\multirow{2}{*}{No. of variables}} & \multicolumn{2}{c|}{\multirow{2}{*}{No. of constraints}} & \multicolumn{2}{c}{Algorithm runtime} \\
			& Arcs & Lower limit \authar[($\underline{\omega}$)] & Upper limit \authar[($\overline{\omega}$)] & \multirow{2}{*}{for TSN} & & & & & & \multicolumn{2}{c}{(secs)}\\
			&  & (minutes) & (minutes) & & & TSN & Proposed\tnote{a} & TSN & Proposed\tnote{b} & TSN & Proposed \\
			\hline
			\multirow{4}{*}{Ex-I\authar[a]} & $a_1$ & 4 & 5 & \multirow{4}{*}{8} & \multirow{4}{*}{1} & \multirow{4}{*}{73} & \multirow{4}{*}{11$_L$+21$_I$} & \multirow{4}{*}{31} & \multirow{4}{*}{35} & \multirow{4}{*}{0.8}  & \multirow{4}{*}{0.7} \\
			& $a_2$ & 4 & 5 & & & & & & & & \\
			& $a_5$ & 7 & 8 & & & & & & & & \\
			& $a_6$ & 7 & 8 & & & &  & & & &\\
			\hline
			\multirow{2}{*}{Ex-I\authar[b]} & $a_4$ & 4 & 5 & \multirow{2}{*}{20} & \multirow{2}{*}{1} & \multirow{2}{*}{188} & \multirow{2}{*}{11$_L$+21$_I$} & \multirow{2}{*}{67} & \multirow{2}{*}{35} & \multirow{2}{*}{0.8}  & \multirow{2}{*}{0.8} \\
			& $a_5$ & 5 & 14 & & & & & & & & \\
			\hline
			\multirow{4}{*}{Ex-I\authar[c]} & $a_3$ & 4.1 & 4.9 & \multirow{4}{*}{12} & \multirow{4}{*}{2} & \multirow{4}{*}{202} & \multirow{4}{*}{22$_L$+42$_I$} & \multirow{4}{*}{84} & \multirow{4}{*}{67} & \multirow{4}{*}{0.8}  & \multirow{4}{*}{0.9} \\
			& $a_3$ & 8.1 & 10.9 & & & & & & & & \\
			& $a_4$ & 4.1 & 4.9 & & & & & & & & \\
			& $a_5$ & 5.5 & 14 & & & & & & & & \\
			\hline
			Ex-II\authar[a] & \multicolumn{3}{c|}{$719$ unavailabilities} & 222 & 2 & 29,532 & \multirow{3}{*}{724$_L$} & 16,069 & \multirow{4}{*}{2901} & 255 & 251 \\
			Ex-II\authar[b] & \multicolumn{3}{c|}{$719$ unavailabilities} & 296 & 2 & 42,434 & \multirow{3}{*}{+ 1084$_I$} & 21,397 & & 2609 & 251 \\
			Ex-II\authar[c] & \multicolumn{3}{c|}{$4727$ unavailabilities} & 1554 & 2 & - & & - & & - & 260 \\
			Ex-II\authar[d] & \multicolumn{3}{c|}{No unavailability} & - & 2 & - & & - & & - & 1.0 \\
			\hline
		\end{tabular}
		\begin{tablenotes}
			\item[a] The subscript $L$ indicates linear (continuous) and $I$ integer.
			\item[b] The constraints due to unavailabilities are not shown these constraints are not used directly (see Algorithm \ref{algo:unav}).
		\end{tablenotes}
		\end{threeparttable}
	\end{center}
    \end{table}
\end{landscape}

The number of decision variables in the TSN representation of the example network is given by the expression $\vert \mathcal{K} \vert \{ (\vert V \vert + \vert A_d \vert + 3) T - 1 - \textit{unav} \}$; where the network contains: $\vert V \vert T$ waiting arcs, $\vert A_d \vert T$ deadhead arcs, $2 T - 3$ service arcs, $T + 1$ sink arcs, and one source arc. $T$ denotes the time limit for TSN, and $\textit{unav}$ denotes the number of arcs that can't be traversed due to unavailability, and therefore doesn't contribute to the number of decisions. Similarly, the number of constraints in TSN is determined by the $\vert \mathcal{K} \vert \{ \vert V \vert (T + 1) + 2\}$ flow constraints at each vertex and $2$ the service arc constraints for the two service arcs; thereby resulting in: $\vert \mathcal{K} \vert \{ \vert V \vert (T + 1) + 2\} + 2$ total constraints. In Figure \ref{fig:tsn3v6a}, there are $6$ unavailable arcs, thereby resulting in ($10 \times 8- 1 - 6 =$) $73$ decision variables, as shown in \textit{Ex-Ia} of Table \ref{tab:table1}. Similarly, there are $29$ flow constraints and $2$ service arc constraints, resulting in a total of $31$ constraints. On the other hand, for the proposed method, there are $\vert K \vert \{ \vert V \vert (\vert A_R \vert + 1) + 2 \}$ replicated graph vertices and $\vert K \vert \{(\vert A_d \vert + 1) (\vert A_R \vert + 1) + \vert A_R \vert^2 + 2\}$ replicated graph arcs. This results in $11$ continuous variables and $21$ integer variables. The number of constraints in the proposed method is the sum of flow constraints for each vertex, running-time constraints for each arc, one temporal constraint for computing maximum completion time, and $2$ service arc constraints. This results in a total of ($11 + 21 + 1 + 2 = $) $34$ constraints. 

Notice that the concept of bound/limit on time is not assumed anywhere in the proposed method, unlike TSN where algorithm convergence is highly sensitive to the number of discrete intervals. This is observed by the change in algorithmic runtime of the two methods as time limit is increased (see Table \ref{tab:table1}, \textit{Ex-IIa} and \textit{Ex-IIb}, KTVK network). The details of the KTVK network is provided in the simulation case study, see Section \ref{rpptru:sec:case}. Note that, time limit for TSN is planner's choice, see $3$rd column of Table \ref{tab:table1}. For simplicity we have chosen a value using prior knowledge of the solution to obtain a significantly small number of time-space decision variables for TSN. 

\subsubsection*{Size comparison over randomly generated graphs} 

The size of any formulation (number of variables and constraints) is one of the principal factors in determining the computation time of NP-Hard problems. We perform computational experiments on the random graph dataset to determine the impact of vertex size on the size of the two formulations: the proposed formulation, and TSN. 
As shown in Figure \ref{fig:tsn_compare}, it is observed that the proposed formulation scales comparatively well with the size of the network in comparison to the TSN based approach. 
Here, the size of the decisions in the TSN representation is approximated here using $(\vert V \vert + \vert A \vert + 1) (\vert A \vert + 1)$; where we assume the $\vert V \vert$ waiting arcs, $\vert A \vert$ regular arcs, and one arc from depot vertex to the end vertex (see TSN representation in Figure \ref{fig:tsn3v6a}) to be repeated $\vert A \vert + 1$ times. The equation results in $70$ decision variables for TSN representation, which is an underestimate to the true value of $73$ shown in Figure \ref{fig:tsn3v6a} and \textit{Ex-Ia} of Table \ref{tab:table1}. The number of constraints for TSN is chosen as $\vert \mathcal{K} \vert \{ \vert V \vert (\vert A \vert + 1) + 2\} + 2$, resulting in $25$ constraints which is another underestimate of the true value of $31$, see \textit{Ex-Ia} in Table \ref{tab:table1}. For the proposed approach, the number of decisions is determined by creating their respective replicated graphs, while the number of constraints (flow, running-time, and service) is two counts larger than the number of decision variables, as discussed in the previous subsection. 

\begin{figure}[!h]
	\centering
	\includegraphics[width=0.98\textwidth]{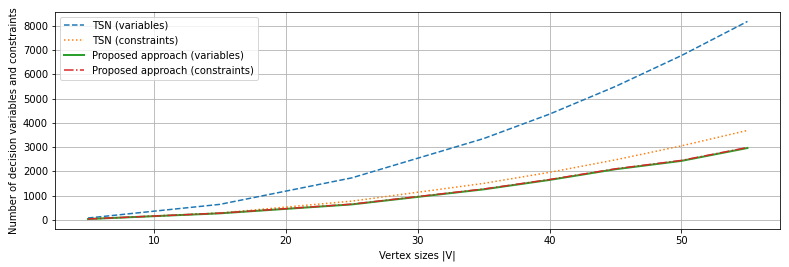}
	\captionsetup{justification=centering}
	\caption{The number of decision variables in TSN grows much faster than the proposed RPP-TRU formulation. For graphs with $5$ vertices, the number of decision variables in TSN is almost twice of the number in RPP-TRU. The number of constraints in TSN also grows faster than the proposed RPP-TRU formulation. }
	\label{fig:tsn_compare}
\end{figure}

\subsection{Comparison study of the two proposed approaches, and with a non-temporal solution} 

The key difference between the two approaches of our proposed algorithm is the introduction of the degree constraints to eliminate all the feasibility-cuts in the second approach. In this subsection, we perform comparison study of the two approaches to analyze their relative performance. Simultaneously, we also present a comparison study between RPP-TRU base-case and RPP-TRU solutions. A comparison study with delayed RPP solution (ignoring the unavailability data) is also demonstrated to show significant cost and time savings by avoiding delay-prone routes. These solutions are denoted as RPP-TRU base-case in the succeeding subsections as it ignores the unavailability restrictions and thereby solves an RPP only. These results are obtained in the first iteration of both the proposed approaches. 

Table \ref{tab:compr} shows all the data related to the comparison study of the two proposed approaches, and the cost savings provided by an optimal/sub-optimal RPP-TRU solution instead of utilizing an RPP based optimal solution. Our results are generated in a HP Probook (Core i7, 8th generation) and 8 GB of memory using Matlab. 

\textbf{The two proposed approaches:} The two approaches of the proposed algorithm are observed to have identical performance parameters over small-medium size graph networks $(\vert V \vert \leq 60)$, see Figures \ref{fig:comp_time} and \ref{fig:niter}. There are differences observed in the computational time in Figure \ref{fig:comp_time}, however there is no clear pattern of one approach dominating the other. In particular, addition of degree constraints to eliminate the family of feasibility-cuts doesn't provide significant improvement in computational efficiency. 

\begin{figure}[!h]
	\centering
	\includegraphics[width=0.98\textwidth]{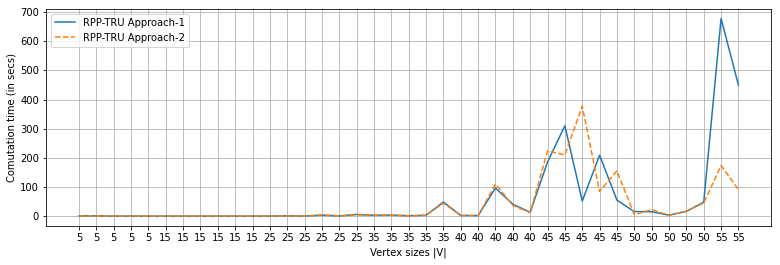}
	\captionsetup{justification=centering}
	\caption{Computation time of the two approaches to generate the RPP-TRU solutions. }
	\label{fig:comp_time}
\end{figure}

\begin{landscape} 
\begin{table}[!htp]
\centering
\caption{Comparing results over randomly generated graphs. }\label{tab:compr}
\scriptsize
\begin{tabular}{lrrrr||c|ccc|ccc||cc|cc|cc}\toprule
$\vert V \vert$ & $\vert A \vert$ & $\vert A_R \vert$ & $\vert \mathcal{V} \vert$ & $\vert \mathcal{A} \vert$ & Cost  & \# iter & Cost & Time & \# iter & Cost & Time & RPP & RPP & A1 & A1 & A2 & A2 \\
& & & & & base-case  & A1 & A1 & A1 & A2 & A2 & A2 & spatial cost & delay & spatial cost & delay & spatial cost & delay \\\midrule
5 &6 &2 &17 &26 &28 &1 &28 &0.32 &1 &28 &0.22 &14 &0 &14 &0 &14 &0 \\
5 &6 &2 &17 &26 &34 &1 &34 &0.26 &1 &34 &0.26 &17 &0 &17 &0 &17 &0 \\
5 &6 &2 &17 &26 &41 &2 &41 &0.28 &2 &41 &0.23 &18 &5 &18 &5 &18 &5 \\
5 &6 &2 &17 &26 &41 &3 &41 &0.26 &3 &41 &0.22 &19 &3 &19 &3 &19 &3 \\
5 &6 &2 &17 &26 &50 &1 &50 &0.21 &1 &50 &0.18 &25 &0 &25 &0 &25 &0 \\
15 &18 &6 &107 &170 &104 &3 &104 &0.34 &3 &104 &0.31 &51 &2 &51 &2 &51 &2 \\
15 &18 &6 &107 &170 &141 &3 &141 &0.38 &3 &141 &0.34 &64 &13 &64 &13 &64 &13 \\
15 &18 &6 &107 &170 &116 &2 &116 &0.32 &2 &116 &0.28 &57 &2 &57 &2 &57 &2 \\
15 &18 &6 &107 &170 &137 &4 &133 &0.38 &4 &133 &0.35 &65 &7 &65 &3 &65 &3 \\
15 &18 &6 &107 &170 &88 &2 &88 &0.30 &2 &88 &0.28 &40 &8 &40 &8 &40 &8 \\
15 &18 &6 &107 &170 &105 &2 &105 &0.31 &2 &105 &0.27 &52 &1 &52 &1 &52 &1 \\
25 &30 &9 &252 &392 &164 &1 &164 &0.44 &1 &164 &0.40 &82 &0 &82 &0 &82 &0 \\
25 &30 &9 &252 &392 &224 &13 &220 &1.06 &1 &220 &0.45 &110 &4 &110 &-0 &110 &-0 \\
25 &30 &9 &252 &392 &260 &2 &246 &0.70 &4 &246 &0.81 &123 &14 &123 &0 &123 &0 \\
25 &30 &9 &252 &392 &191 &41 &184 &6.09 &41 &184 &5.68 &89 &13 &89 &6 &89 &6 \\
25 &30 &9 &252 &392 &212 &3 &202 &0.91 &4 &202 &0.90 &98 &16 &101 &0 &101 &0 \\
25 &30 &9 &252 &392 &194 &22 &194 &6.04 &22 &194 &5.31 &94 &6 &94 &6 &94 &6 \\
35 &42 &13 &492 &772 &251 &7 &251 &2.85 &7 &251 &3.58 &120 &11 &120 &11 &120 &11 \\
35 &42 &13 &492 &772 &251 &5 &239 &3.63 &5 &239 &2.95 &108 &35 &108 &23 &108 &23 \\
35 &42 &13 &492 &772 &288 &4 &282 &1.40 &4 &282 &1.33 &138 &12 &138 &6 &138 &6 \\
35 &42 &13 &492 &772 &352 &9 &352 &2.84 &9 &352 &3.01 &167 &18 &167 &18 &167 &18 \\
35 &42 &13 &492 &772 &245 &62 &245 &50.89 &62 &245 &52.38 &116 &13 &116 &13 &116 &13 \\
40 &48 &15 &642 &1010 &289 &4 &287 &2.06 &4 &287 &2.05 &143 &3 &143 &1 &143 &1 \\
40 &48 &15 &642 &1010 &295 &3 &268 &1.64 &3 &268 &1.53 &133 &29 &133 &2 &133 &2 \\
40 &48 &15 &642 &1010 &302 &100 &302 &119.29 &100 &302 &124.39 &141 &20 &141 &20 &141 &20 \\
40 &48 &15 &642 &1010 &415 &38 &378 &41.04 &38 &378 &37.65 &186 &43 &186 &6 &186 &6 \\
40 &48 &15 &642 &1010 &279 &11 &265 &12.64 &11 &265 &12.17 &124 &31 &124 &17 &124 &17 \\
45 &54 &17 &812 &1280 &380 &44 &364 &51.63 &44 &364 &377.09 &176 &28 &176 &12 &176 &12 \\
45 &54 &17 &812 &1280 &393 &27 &375 &209.20 &27 &375 &84.67 &174 &45 &174 &46 &174 &46 \\
45 &54 &17 &812 &1280 &437 &27 &410 &54.92 &52 &410 &155.64 &205 &27 &205 &-0 &205 &-0 \\
50 &60 &18 &952 &1484 &418 &22 &404 &15.74 &9 &404 &5.23 &202 &14 &202 &-0 &202 &-0 \\
50 &60 &18 &952 &1484 &313 &5 &313 &15.18 &5 &313 &22.53 &154 &5 &154 &5 &154 &5 \\
50 &60 &18 &952 &1484 &391 &5 &368 &3.58 &2 &368 &2.43 &184 &23 &184 &0 &184 &0 \\
50 &60 &18 &952 &1484 &407 &11 &403 &16.60 &11 &403 &17.35 &195 &17 &197 &23 &197 &23 \\
50 &60 &18 &952 &1484 &361 &13 &361 &47.85 &13 &361 &44.83 &168 &25 &168 &25 &168 &25 \\
55 &66 &20 &1157 &1808 &445 &100 &445 &678.12 &100 &445 &173.81 &210 &25 &212 &35 &212 &35 \\
55 &66 &20 &1157 &1808 &406 &66 &386 &449.84 &65 &386 &91.38 &187 &32 &187 &22 &187 &22 \\
\bottomrule
\end{tabular}
\end{table}
\end{landscape}

\begin{figure}[!h]
	\centering
	\includegraphics[width=0.98\textwidth]{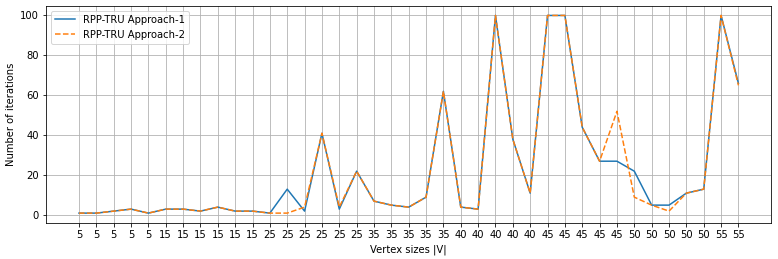}
	\captionsetup{justification=centering}
	\caption{Number of iterations (at most 100) required to generate the solution. }
	\label{fig:niter}
\end{figure}

\textbf{With RPP-TRU base-case:} The RPP-TRU base-case depicts non-temporal RPP solution with a low cost optimal route, and it doesn't account for the delay caused by unavailabilities. Interestingly, in most of the case for small and moderate size medium-dense graphs, it is observed that there is an improvement in the time taken to complete the task at no increase in the spatial/fuel cost, see Figure \ref{fig:spat_only}. In particular, the solution of RPP-TRU has same cost as that of RPP, however a lot of travel time is reduced due to less delay, see Figure \ref{fig:delay_only}. Figure \ref{fig:comp_time} shows the trade-off observed in the computational time to generate these time-saving plans. The percentage improvement in the spatio-temporal costs of the two approaches are shown in Figure \ref{fig:sptmp_cst}. 

\begin{figure}[!h]
	\centering
	\includegraphics[width=0.98\textwidth]{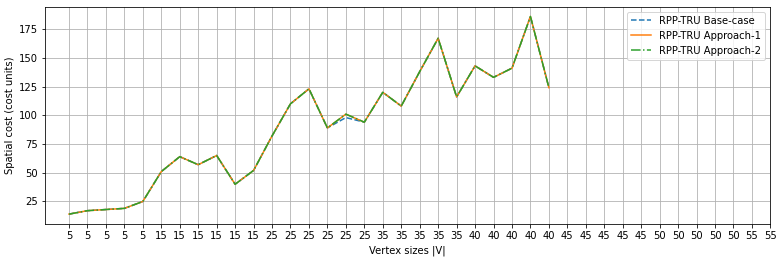}
	\captionsetup{justification=centering}
	\caption{Comparison of the results showing spatial cost only for the solutions from RPP-TRU base-case and proposed RPP-TRU approaches. }
	\label{fig:spat_only}
\end{figure}
\begin{figure}[!h]
	\centering
	\includegraphics[width=0.98\textwidth]{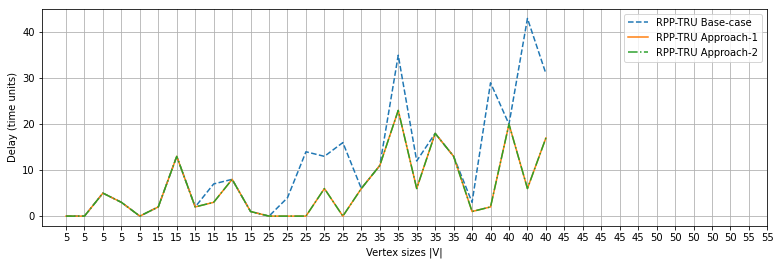}
	\captionsetup{justification=centering}
	\caption{Comparison of the results showing delay or waiting-time observed in RPP-TRU base-case solutions and RPP-TRU solutions.  }
	\label{fig:delay_only}
\end{figure}
\begin{figure}[!h]
	\centering
	\includegraphics[width=0.98\textwidth]{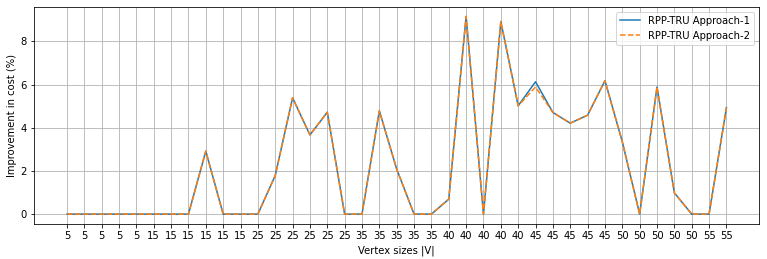}
	\captionsetup{justification=centering}
	\caption{Percentage improvement in the spatio-temporal cost for network of various sizes w.r.t. the RPP-TRU base-case solution. }
	\label{fig:sptmp_cst}
\end{figure}


\subsection{Simulation case study: KTVK railway network}
\label{rpptru:sec:case}

Kurla-Thane-Vashi-Kurla (KTVK) sub-urban railway network is located in Mumbai, India. Figure \ref{fig:ktvk} shows the network map as well as a graph representation with $36$ vertices and $45$ arcs. Among them $9$ arcs, highlighted with bold arrows, model the railway track sections that require servicing (inspection). These service arcs have parallel deadhead arc representations to model motion of the agents on these tracks without performing any servicing; resulting in a total of $54$ arcs. Four types of train movements are assumed: Kurla-Thane fast and slow to-fro (between $v_1$ and $v_3$), Kurla-Vashi slow to-fro (between $v_1$ and $v_4$) and Thane-Vashi slow to-fro (between $v_1$ and $v_4$). The train schedules are assumed to be repeating every $74$ minutes. 

\begin{figure*}[!h]
	\centering
	\captionsetup{justification=centering}
	\includegraphics[width=0.88\textwidth]{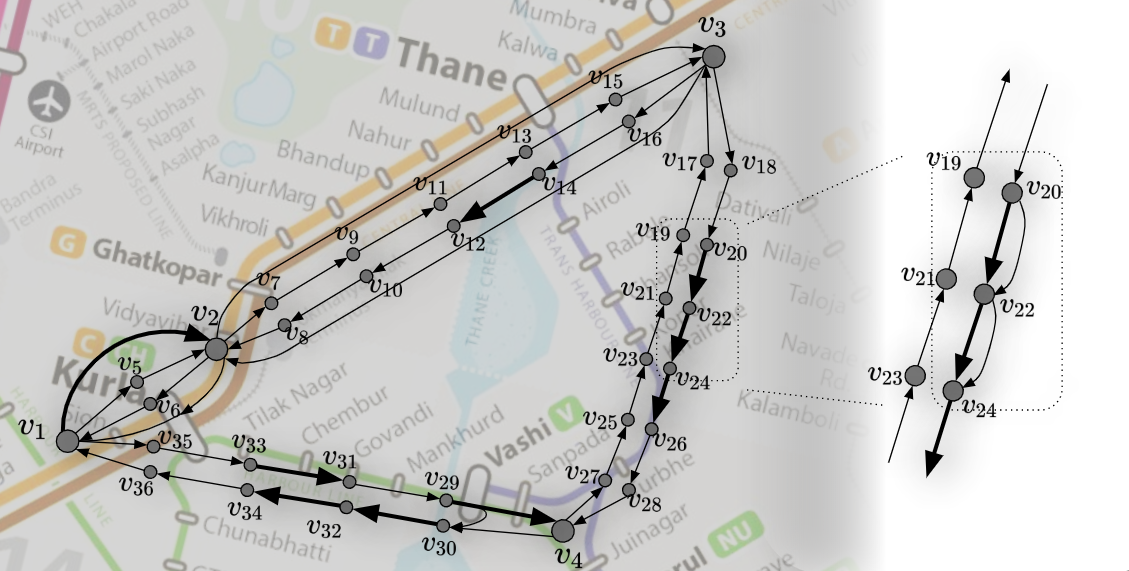}
	\caption{This graph portrays 36 stations all modeled as vertices - indicated using circles of two sizes. There are 45 deadhead arcs in this model, of which 9 are to be serviced. As shown in the enlarged figure in right, the service arcs marked in bold have a duplicate deadhead arc in parallel.}
	\label{fig:ktvk} 
\end{figure*}

A directed multi-graph is constructed to resemble the KTVK suburban railway network. Since the entire network layout data is not available with us, following assumptions are made: all stations considered in this network can park any number of agent for indefinite period of time, unavailability schedules are periodic, and its possible for an agent to move between any two stations (while servicing or cruising) in any $74$ minute time interval i.e. $T_p$-recurrence is guaranteed for the data; where $T_p$ is $74$ minutes. Such periodic schedule that repeats every $74$ minutes is not actually employed by Mumbai sub-urban railways, however such schedules are implemented in different places around the globe - e.g. Dutch railways has timetable that repeats every hour \cite{peeters}. 

\textit{Ex-II} in Table \ref{tab:table1} shows the simulation results over the KTVK network for a variety of unavailability settings. In \textit{Ex-IIa}, a scenario is considered with prior knowledge of the total inspection time taken by optimal solution (207 minutes). Hence the least integral multiple of $T_p$ that can contain the optimal solution, i.e. 3 times $T_p$ ( $= 222$ minutes), is considered as time limit for formulating TSN. In \textit{Ex-IIb}, effect on increasing the total inspection time is observed, by considering 4 times $T_p$ ($= 296$ minutes) as the time-limit for formulating TSN. Using Corollary \ref{cor:timebound}, an inspection time bound of $(\vert V \vert \vert \mathcal{L} \vert + 1) \times (T_p + W_{max}) = 30324$ is obtained; where $\vert V \vert = 36$, $\vert \mathcal{L} \vert = \vert A_R \vert + 1 = 10$, $T_p = 74$, and $W_{max} = 10$. In \textit{Ex-IIc}, \textit{out of memory} error was obtained for the time bound of $1554$, which is much less than $30324$. Lastly, in \textit{Ex-IId}, simulation without any unavailability scenario is presented. 
For the special case (\textit{Ex-IId}), there are existing RPP algorithms to provide a comparable, if not better, solution. In practice, since the two agents are identical, the trajectories of the two agents can be interchanged if there is an overtaking occurrence in the schedule. Also note that, the comparison study with TSN shows that unavailability has a significantly small impact on the proposed method. 

\ifsubbuild

\bibliographystyle{cas-model2-names}

\bibliography{citationsall2.bib}

\fi

\section*{Conclusion} 
\label{rpptu:sec:concl}

In this work, we show that routing and scheduling over network with time-dependent unavailability scenarios can be planned optimally using the proposed methodology. In particular, for the general setting of multi-agent temporal RPP, our work contributes towards methodology for obtaining exact solutions that doesn't disrupt the existing schedules that are modeled as unavailability constraints. Consequently, in a railway application scenario, the resulting inspection schedules get sandwiched between the regular train schedules, thereby allowing frequent monitoring of railway tracks. This improves the safety and reliability of the railway services. 
The simulation case study, presented in this work, demonstrates such an application over KTVK railway network. 

The backbone of this work is established by Theorem \ref{thm:wd} which shows the existence of an optimal solution in the replicated graph-based formulation. In particular, the assumption of well-defined and $T_p$-recurrent graph-agent data dictates the ability to comprehensively model and solve networks with temporal attributes. By insisting the temporal graph-agent data to be well-defined, compatibility is enforced between the network and the agents i.e. it eliminates scenarios with slow agents on a busy network. In particular, an agent that moves too slowly (e.g. ultrasonic inspection vehicles) will not be able to traverse most tracks of a sub-urban railway network during peak hours. Having well-defined temporal graph-agent data, or producing one before attempting the optimization problem, is a wise choice to reduce the computational burden. Many applications have a periodic component in their existing schedules, hence $T_p$-recurrence is an intrinsic property in these applications e.g. railways. 
The comparison results with the TSN approach demonstrate the benefit of the proposed replicated graph-based formulation with continuous time variables over discretized time-space models. Unlike TSN, the proposed methodology is observed to be only slightly affected by the size of the unavailability schedules. 

Two approaches to the proposed algorithm are developed, implemented, and analyzed in this work. The algorithm is developed from Benders' decomposition framework, thereby producing an exact solution given sufficient computational time and memory. In case of an early termination of the algorithm, sub-optimal solutions are obtained along with a dualtiy gap to certify/quantify the level of sub-optimality. Thus the proposed methodology is suited to on-demand routing and scheduling requirements, where the proposed methodology produces either optimal solutions or sub-optimal solutions with duality-bound guarantees while ensuring no disruption of the existing schedules. The results show that the two approaches have similar computational performance over medium-sized sparse networks. In comparison with RPP solutions that ignore existing schedules during the planning phase, the proposed approaches produce good cost and time improvement at the cost of computational time. 
In particular, there is negligible compromise in the cost of traversal, while producing moderately good improvement in the planned schedules (less delay/waiting). The trade-off is observed in the computational time of the algorithms, which may be compensated by solving medium-sized dense networks in a relatively powerful computational system. 

The proposed work has real-world applications over macroscopic models of railway networks. Mesoscopic models pose additional challenges arising from practical side constraints. Implementation-oriented study and analysis over mesoscopic model to incorporate constraints like speed limit, headway, etc is one useful research direction for the future. A study directed towards robust schedules i.e. minimum impact of delays on the planned schedules may be another useful research direction that addresses the challenge of schedule uncertainty.
In this work, an indefinite waiting permit is assumed at the vertices which restricts our result to cases with a small number of agents in many applications. Since a large agent swarm is not practical in railways, and considering the fact that railway network problems are solved regionally by railway undertakings, portions of railway networks can be modeled with small and sparser graphs, making our results suitable for regional implementation. Swarm cases on large and dense graphs can be considered a future work to inspect complex railway networks on a large scale. 


\printbibliography[title=References]

\end{document}